\documentclass[11pt]{amsart}

\usepackage{amsfonts,amssymb,amsmath,mathrsfs,graphicx}
\usepackage{float}
\usepackage{epsfig}

\usepackage{CJK,color}

\definecolor{red}{rgb}{1.00,0.00,0.00}
\definecolor{blue}{rgb}{0.00,0.00,0.63}
\definecolor{black}{rgb}{0.00,0.00,0.00}

\newtheorem{theorem}{Theorem}[section]
\newtheorem{lemma}{Lemma}[section]

\newtheorem{proposition}{Proposition}[section]
\newtheorem{remark}{Remark}[section]

\newcommand{\ba}{\begin{aligned}}
\newcommand{\ea}{\end{aligned}}
\newcommand{\be}{\begin{equation}}
\newcommand{\ee}{\end{equation}}
\renewcommand{\div}{ {\rm div }  }
\def\na{\nabla}

\def\charf {\mbox{{\text 1}\kern-.30em {\text l}}}

\def\lam{\lambda}  

\def\f{\frac}
\def\t{\theta}
\def\d{\delta}

\def\z{\zeta}

\def\di{\displaystyle}
\def\r{\rho}

\def\k{\kappa}

\def\v{\varepsilon}

\def\charf {\mbox{{\text 1}\kern-.30em {\text l}}}

\def\lam{\lambda}  

\def\f{\frac}
\def\t{\theta}

\def\d{\delta}
\def\s{\sigma}

\def\di{\displaystyle}
\def\r{\rho}

\def\k{\kappa}
\def\pa{\partial}

\begin{document}

\title[Stability of planar rarefaction wave to CNS equations]
{Stability of planar rarefaction wave to 3D full compressible Navier-Stokes equations}

\author[Li]{Lin-an Li}
\address[Lin-an Li]{\newline Institute of Applied Mathematics,AMSS, Chinese Academy of Sciences, Beijing 100190, P. R. China
\newline and School of Mathematical Sciences, University of Chinese Academy of Sciences, Beijing 100049, P. R. China}
\email{linanli@amss.ac.cn}

\author[Wang]{Teng Wang}
\address[Teng Wang]{\newline College of Applied Sciences, Beijing University of Technology, Beijing 100124, P. R. China}
\email{tengwang@amss.ac.cn}

\author[Wang]{Yi Wang}
\address[Yi Wang]{\newline  CEMS, HCMS, NCMIS, Academy of Mathematics and Systems Science, Chinese Academy of Sciences, Beijing 100190, P. R. China
\newline and School of Mathematical Sciences, University of Chinese Academy of Sciences, Beijing 100049, P. R. China
}
\email{wangyi@amss.ac.cn}

\date{\today}

\maketitle

\begin{abstract}

We prove the time-asymptotic stability toward planar rarefaction wave for the
three-dimensional full compressible Navier-Stokes equations in an infinite long flat nozzle domain $\mathbb{R}\times\mathbb{T}^2$.
Compared with one-dimensional case,  the proof here is based on our new observations on the cancellations on the flux terms and viscous terms due to the underlying wave structures, which are
crucial to overcome the difficulties due to the wave propagation along the transverse directions $x_2$ and
$x_3$ and its interactions with the planar rarefaction wave in $x_1$ direction.

\end{abstract}

\maketitle \centerline{\date}


%
%
\section{Introduction}
\setcounter{equation}{0}
The motion of compressible viscous and heat-conductive  fluid occupying a spatial domain $\Omega\subset \mathbb{R}^3$
is governed by the following full compressible Navier-Stokes-Fourier system:
\begin{equation}\label{ns}
\begin{cases}
\di \r_t+\div(\r u)=0,\\
\di (\r u)_t+\div(\r u\otimes u)+\na p=\div\mathcal{T},\\
\di (\r E)_t+\div (\r E u + p u)=\kappa\Delta\t+\div(u\mathcal{T}),
\end{cases}
\end{equation}
where $t\geq0$ is the time variable and $x=(x_1,x_2,x_3)\in \Omega$ is the spatial variables. In the present paper, we are concerned with the viscous fluid flowing in an infinite long flat nozzle domain $\Omega:=\mathbb{R}\times\mathbb{T}^2$ with $\mathbb{R}$ being a real line and
$\mathbb{T}^2:= (\mathbb{R}/ \mathbb{Z})^2$ being a two-dimensional unit flat torus. The functions $\r, u=(u_1,u_2,u_3)^t, p$ and $\t$ represent respectively the fluid density, velocity, pressure and absolute temperature and $E=e+\f12|u|^2$ is the specific total energy with $e$ being the internal energy, and $\mathcal{T}$ is the viscous stress tensor given by
\be
\mathcal{T}= 2\mu \mathbb{D}(u) + \lambda\div u \mathbb{I}
\ee
where $\mathbb{D}(u)=\frac{\na u + (\na u)^t}{2}$ stands for the deformation tensor, $\mathbb{I}$ is the $3\times3$ identity matrix and $\mu$ and $\lam$ represent the shear and bulk viscosity coefficients of the fluids respectively and they both are constants satisfying the physical constraints:
\be
\mu>0,\quad 2\mu+3\lambda\geq0.
\ee
Moreover, the constant $\kappa>0$ denotes the heat-conductivity coefficient for the fluids. The equations \eqref{ns} then express respectively the conservation of mass, the balance of momentum, and the balance of energy for the flow under the effect of the inner pressure, viscosities and the conduction of thermal energy. Here we investigate the ideal poly-tropic fluids such that the pressure $p$ and the internal energy $e$ are given by the following state equations:
\be
p=R\r\t=A\r^{\gamma}\exp\Big(\f{\gamma-1}{R}S\Big),\quad e=\f{R}{\gamma-1}\t,
\ee
where $S$ is the entropy, $\gamma>1$ is the adiabatic exponent, and both $A$ and $R$ are positive fluid constants.

\
The following initial data is imposed to the system \eqref{ns}
\begin{equation}\label{initial}
(\r,u,\t)(x,t=0)=(\r_0,u_0,\t_0)(x),\qquad x\in \Omega.
\end{equation}
Since we are concerned with the stability of the planar rarefaction wave to the system \eqref{ns}, we consider the following far-fields conditions on the $x_1$-direction
\be\label{far}
(\r,u,\t)(x,t) \to (\r_\pm,u_\pm,\t_\pm), \quad {\rm as}\quad x_1\rightarrow \pm\infty,~ t>0,
\ee
with $u_\pm=(u_{1\pm},0,0)^t$ and $\rho_\pm>0, ~u_{1\pm}, ~\t_\pm>0$ are prescribed constant states,
and the periodic boundary conditions are imposed on $(x_2,x_3)\in\mathbb{T}^2$ for the solution $(\r,u,\t)(x,t)$. Moreover, the two end states $(\rho_\pm,u_{\pm},\t_\pm)$ are connected by the rarefaction wave solution to the Riemann problem of the corresponding 1D compressible Euler system:
\begin{equation}\label{euler}
\left\{
\begin{array}{l}
\di\rho_t+(\rho u_1)_{x_1}=0,\\
\di(\rho u_1)_t+(\rho u_1^2+p)_{x_1}=0,\\
\di (\rho E)_t+(\rho Eu_1+pu_1)_{x_1}=0,
\end{array}
\right.
\end{equation}
with the Riemann initial data
\begin{equation}\label{R-in}
(\r,u_1,\t)(x_1,0)=(\r_0^r,u_{10}^r,\t_0^r)(x_1)=\left\{
\begin{array}{ll}
\di (\r_-, u_{1-}, \t_-), &\di x_1<0, \\
\di (\r_+, u_{1+}, \t_+), &\di x_1>0.
\end{array}
 \right.
\end{equation}

\

It could be expected that the large-time behavior of the solution to the compressible Navier-Stokes equations
\eqref{ns}-\eqref{far} is closely related to the Riemann problem to the corresponding three-dimensional compressible Euler equations
\begin{equation}\label{euler'}
\left\{
\begin{array}{l}
\di\rho_t+\div(\rho u)=0,\\
\di(\rho u)_t+\div(\rho u\otimes u)+\na p=0,\\
\di(\r E)_t+\div (\r Eu+ p u)=0,
\end{array}
\right.
\end{equation}
with the Riemann initial data
\begin{equation}\label{R-in'}
(\r,u,\t)(x,0)=(\r_0^r,u_{0}^r,\t_0^r)(x)=\left\{
\begin{array}{ll}
\di (\r_-, u_{-}, \t_-), &\di x_1<0 \\
\di (\r_+, u_{+}, \t_+), &\di x_1>0.
\end{array}
 \right.
\end{equation}

\

The inviscid compressible Euler system \eqref{euler} or \eqref{euler'} is an ideal fluid model with the dissipative effects being neglected, which is a typical example of the system
of hyperbolic  conservation laws.  The most important feature of  the hyperbolic system \eqref{euler} or \eqref{euler'}  is that its classical solution may blow up, that is, the shock may form, in finite time, no matter how smooth or small the initial data is. Actually, there are three basic wave patterns to the system of hyperbolic conservation laws,
i.e., shock and rarefaction waves in the genuinely nonlinear characteristic fields, and contact discontinuity in the linearly degenerate fields.
However, the motion of real fluids should take into account the effects of both viscosities and
heat-conductivity, which is described by the compressible Navier-Stokes system \eqref{ns}, i. e., the corresponding viscous system of inviscid Euler system \eqref{euler} or \eqref{euler'}.
Moreover, it can be expected that the large-time behavior of the solutions to the compressible Navier-Stokes equations \eqref{ns}-\eqref{far} is governed by the
solutions of the corresponding Riemann problem \eqref{euler'}-\eqref{R-in'}, which contains planar shock wave, planar rarefaction wave and
contact discontinuity in general.  It is interesting and important to investigate the time-asymptotic stability
of these basic planar wave patterns to the compressible Navier-Stokes equations \eqref{ns} in higher dimension. In the present paper, we first study the nonlinear stability of planar rarefaction wave to the system \eqref{ns} in an infinite long flat nozzle domain $\Omega:=\mathbb{R}\times\mathbb{T}^2$.

\

On one hand, there are essential differences between the one-dimensional Riemann problem \eqref{euler}-\eqref{R-in} and the multi-dimensional Riemann problem \eqref{euler'}-\eqref{R-in'} even with the components $u_2$, $u_3$ are continuous on both sides of $x_1=0$ as in \eqref{R-in'}. Precisely speaking, it is first proved by Chiodaroli, De Lellis and Kreml \cite{CDK-1} and Chiodaroli and Kreml \cite{CK} that there exist infinitely many bounded admissible weak solutions to \eqref{euler'}-\eqref{R-in'} in two-dimensional isentropic regime satisfying the natural entropy condition for shock Riemann initial data by using the elegant convex integration methods as in De Lellis and Sz$\acute{{\rm e}}$kelyhidi \cite{DS}.  Meanwhile, the construction of weak solutions in \cite{CDK-1,CK} seems essential to the two-dimensional system and can not be applied to one-dimensional problem \eqref{euler}-\eqref{R-in}. Then Klingenberg and Markfelder \cite{KM} and Brezina, Chiodaroli and Kreml \cite{BCK} extend the results in \cite{CDK-1,CK} to the case when the corresponding Riemann initial data contain shock or contact discontinuity. On the other hand, similar to the one-dimensional case, for the Riemann solution only containing rarefaction waves to \eqref{euler'}-\eqref{R-in'}, Chen and Chen \cite{CC} and Feireisl and Kreml \cite{FK}, Feireisl, Kreml and Vasseur \cite{FKV} independently proved rarefaction wave is unique in the class of bounded weak solution to \eqref{euler'}-\eqref{R-in'} even the rarefaction waves are connected with vacuum states (cf. \cite{CC}).

\

As mentioned before, the inviscid Euler system \eqref{euler} or \eqref{euler'} is an ideal fluid model and the real fluids could be described by the viscous system  \eqref{ns}, which is a typical example of  the system of the viscous conservation laws.  Deep investigations have been achieved on the nonlinear stability of basic wave patterns for viscous conservation laws in one-dimensional case. For the asymptotic stability of viscous shock profile,
it started from Goodman \cite{goodman} for the uniformly viscous conservation laws and Matsumura and Nishihara \cite{MN-85} for the compressible Navier-Stokes equations with physical viscosities independently by the anti-derivative methods under zero mass condition imposed on the initial perturbation. Then Liu \cite{Liu-1985} and Szepessy and Xin \cite{Szepessy-Xin} removed the zero mass condition for the uniformly viscous conservation laws by introducing  the suitable shift on the shock profile and diffusion waves in the transverse characteristic fields  and Liu and Zeng \cite{Liu-Zeng} for the physical viscosity case. For the stability of rarefaction wave, we refer to Matsumura and Nishihara \cite{MN-86, MN-92} for isentropic compressible Navier-Stokes equations and Liu and Xin \cite{liu-xin-1}, Nishihara, Yang and Zhao \cite{NYZ} for non-isentropic system. Then Liu and Yu \cite{liu-yu-rare} proved the stability of rarefaction wave to one-dimensional general $n\times n$ conservation laws system with artificial viscosity by point-wise Green function methods.
For the stability of viscous contact discontinuity wave, we refer to Liu and Xin \cite{Liu-Xin} and Xin \cite{xin-2} for the uniformly viscous conservation laws and Huang, Matsumura and Xin \cite{Huang-Matsumura-Xin} for the compressible Navier-Stokes equations under the zero mass condition on the perturbation,
Then Huang, Xin and Yang \cite{Huang-Xin-Yang} removed this zero mass condition in \cite{Huang-Matsumura-Xin} for the 1D compressible Navier-Stokes equations \eqref{ns}.
For the composite waves, Huang and Matsumura \cite{Huang-matsumura} first studied the asymptotic stability of two viscous shock waves under
general initial perturbation without zero mass conditions on initial perturbations for the full  1D Navier-Stokes system and
Huang, Li and Matsumura \cite{Huang-Li-matsumura} justified the stability of a combination wave of a viscous contact wave and rarefaction waves.
Recently, Huang and Wang \cite{huang-wang} improved the stability result in \cite{Huang-Li-matsumura} to a class of large initial perturbations.

\

Although there have been rather satisfactory results about the stability of basic wave patterns for viscous conservation laws in the one-dimensional case, the stability toward the planar wave patterns for the compressible Navier-Stokes equations \eqref{ns} in multi-dimensional case is still open due to the higher dimensionality. For the scalar viscous conservation laws, Xin \cite{Xin} proved the asymptotic stability of planar rarefaction waves in multi-dimensional case by elementary $L^2$-energy
method in 1990. Then Ito \cite{Ito} and Nishikawa and Nishihara \cite{NN} extended the stability result in \cite{Xin} by obtaining the decay rate in time. For an artificial $2\times 2$ system with positively definite viscosity matrix, Hokari and Matsumura \cite{HM} proved the stability of the planar rarefaction wave in two-dimensional case, which crucially depends on the strict positivity of the viscosity matrix and can not be applied to the compressible Navier-Stokes system \eqref{ns} with physical viscosities.
For the compressible and isentropic Navier-Stokes equations, which is the special case of the system \eqref{ns} with the entropy being constant and the energy equation can be decoupled and neglected, the first and third author of the present paper Li and Wang \cite{LW} proved the stability of planar rarefaction wave in two-dimensional domain $\mathbb{R}\times\mathbb{T}$.

\

In the present paper, we shall prove the time-asymptotic stability of the planar rarefaction wave for the three-dimensional full compressible Navier-Stokes equation \eqref{ns} with physical viscosities and heat-conductivity for any adiabatic exponent $\gamma>1$. Compared with the one-dimensional stability results
in \cite{MN-86,MN-92,NYZ}, the main difference here lies in higher dimensionality and the physical viscosities terms coupled in momentum equation \eqref{ns}$_2$ and the energy equation \eqref{ns}$_3$  and we can not use the technique for one-dimensional fluid by substituting the mass equation \eqref{ns}$_1$ into the momentum equation \eqref{ns}$_2$ directly to obtain the derivative estimates of the density function as in \cite{MN-86,MN-92,NYZ}. Compared with the two-dimensional stability result for isentropic flow in \cite{LW}, the full compressible Navier-Stokes equation \eqref{ns} here is a real physical model involving the thermal conduction and the main difference lies in the thermal energy equation \eqref{ns}$_3$ additionally in three-dimensional domain.
Fortunately, we observe some cancellations
between the flux terms and viscosity terms for the full compressible Navier-Stokes equations \eqref{ns}  such that we can successfully
overcome the difficulties due to the planar rarefaction wave propagation in $x_2, x_3$-directions and its interactions with $x_1$-direction and finally we can prove our time-asymptotic stability toward the planar rarefaction wave. More precisely, we prove that if the initial data $(\r_0,u_0,\t_0)$ in \eqref{initial} is suitably close to the planar rarefaction wave, then the three-dimensional problem \eqref{ns}-\eqref{far} admits a  global-in-time smooth solution which tends to the planar rarefaction wave as $t\rightarrow+\infty$. Note that the rarefaction wave strength $|(\r_+-\r_-,u_+-u_-,\t_+-\t_-)|$ here need not to be sufficiently small.  The detailed stability result can be found in Theorem \ref{thm} below.

\

To state our main result, we first recall the planar rarefaction wave. It is straight to calculate that the Euler system \eqref{euler} for $(\r,u_1,\t)$ has three distinct eigenvalues
$$
\lambda_i(\r,u_1,S)=u_1+(-1)^{\f{i+1}{2}}\sqrt{p_\r(\r,S)},~ i=1,3,\qquad
\lambda_2(\r,u_1,S)=u_1,~
$$
with corresponding right eigenvectors
$$
r_i(\r,u_1,S)=((-1)^{\f{i+1}{2}}\r,\sqrt{p_\r(\r,S)},0)^t,~i=1,3, \qquad\quad r_2(\r,u_1,S)=(p_{\scriptscriptstyle S},0,-p_\rho)^t,$$
such that
$$
r_i(\r,u_1,S)\cdot \nabla_{(\r,u_1,S)}\lambda_i(\r,u_1,S)\neq 0,~ i=1,3,\quad
{\rm
and}\quad
r_2(\r,u_1,S)\cdot \nabla_{(\r,u_1,S)}\lambda_2(\r,u_1,S)\equiv 0.
$$
Thus the  two $i$-Riemann invariants $\Sigma_i^{(j)}(i=1,3, j=1,2)$ can be defined by (cf. \cite{Smoller})
\begin{equation}\label{RI}
\Sigma_i^{(1)}=u_1+(-1)^{\f{i-1}{2}}\int^{\r}\f{\sqrt{p_z(z,S)}}{z}dz,\qquad
\Sigma_i^{(2)}=S,
\end{equation}
such that
$$
\nabla_{(\r,u_1,S)} \Sigma_i^{(j)}(\r,u_1,S)\cdot r_i(\r,u_1,S)\equiv0,\quad i=1,3,~ j=1,2.
$$
Given the right state $(\rho_+, u_{1+}, \theta_+)$ with $\r_+>0,\t_+>0$,  the $i$-rarefaction wave curve $(i=1,3)$ in the phase space $(\r,u_1,\t)$ with $\r>0$ and $\t>0$ can be defined by (cf. \cite{Lax}):
\begin{equation}\label{R3-curve}
 R_i (\rho_+, u_{1+}, \theta_+):=\Bigg{ \{} (\rho, u_1, \theta)\Bigg{ |}\lambda_{ix_1}(\r,u_1,S)>0, \Sigma_i^{(j)}(\r,u_1,S)=\Sigma_i^{(j)}(\r_+,u_{1+},S_+),~~j=1,2\Bigg{ \}}.
\end{equation}
Without loss of generality, we consider the stability of planar $3-$rarefaction wave to the Euler system \eqref{euler}, \eqref{R-in} in the present paper and the stability of $1-$rarefaction wave can be done similarly.
The $3-$rarefaction wave to the Euler system \eqref{euler}, \eqref{R-in} can be expressed explicitly by
 the Riemann solution to the inviscid Burgers equation:
\begin{equation}\label{bur}
\left\{\begin{array}{ll}
w_t+ww_{x_1}=0,\\
w(x_1,0)=w_0^r(x_1)=\left\{\begin{array}{ll}
w_-,&x_1<0,\\
w_+,&x_1>0.
\end{array}
\right.
\end{array}
\right.
\end{equation}
If $w_-<w_+$, then the Riemann problem $(\ref {bur})$ admits a
rarefaction wave solution $w^r(x_1, t) = w^r(\f {x_1}{t})$ given by
\begin{equation}\label{abur}
w^r\Big(\f {x_1}{t}\Big)=\left\{\begin{array}{lr}
w_-,&\f {x_1}{t}\leq w_-,\\
\f {x_1}{t},&w_-\leq \f {x_1}{t}\leq w_+,\\
 w_+,&\f {x_1}{t}\geq w_+.
\end{array}
\right.
\end{equation}
Then the 3-rarefaction wave solution $(\r^r,u_1^r,\t^r)(\f {x_1}{t})$ to the
compressible Euler equations \eqref{euler}, \eqref{R-in} can be defined explicitly by
\begin{eqnarray}   \label{3-rw}
\left\{
\begin{array}{l}
\di w_\pm=\lambda_3(\r_\pm,u_{1\pm},\t_\pm), \qquad w^r(\f {x_1}{t})= \lambda_3(\r^r,u_1^r,\t^r)(\f {x_1}{t}),\\
\di
\Sigma_3^{(j)}(\r^r,u_1^r,\t^r)(\f{x_1}{t})=\Sigma_3^{(j)}(\rho_\pm,u_{1\pm},\t_\pm),\quad j=1,2,\quad
u_2^r= u_3^r=0,
\end{array} \right.
\end{eqnarray}
where $\Sigma_3^{(j)}~(j=1,2)$ are the 3-Riemann invariants defined in \eqref{RI}.

We construct a smooth 3-rarefaction wave profile to the wave fan defined in \eqref{3-rw}.
Motivated by \cite{MN-92},  the smooth rarefaction wave can be constructed by the Burgers equation
\begin{equation}\label{dbur}
\left\{
\begin{array}{l}
\di \bar w_{t}+\bar w\bar w_{x_1}=0,\\
\di \bar w(x_1,0)=\bar w_0(x_1)=\f{w_++w_-}{2}+\f{w_+-w_-}{2}k_q\int_0^{\varepsilon x_1}(1+y^2)^{-q}dy,
\end{array}
\right.
\end{equation}
where $\v>0$ is a small constant to be determined and $k_q$ is a positive constant such that $k_q\int_0^{\infty}(1+y^2)^{-q}dy=1$ for each $q\geq2$.
Note that the solution $\bar w(x_1,t)$ of the problem \eqref{dbur} can be given explicitly by
\begin{equation}\label{b-s}
\bar w(x_1,t)=\bar w_0(x_0(x_1,t)),\qquad x_1=x_0(x_1,t)+\bar w_0(x_0(x_1,t))t.
\end{equation}

Correspondingly, the smooth rarefaction wave profile $(\bar{\r},\bar{u},\bar{\t})(x_1,t)$ to compressible Euler equations
$\eqref {euler}, \eqref{R-in}$ can be defined by
\begin{eqnarray}
\left\{
\begin{array}{l}
\di w_\pm=\lambda_3(\r_\pm,u_{1\pm},\t_\pm), \qquad \bar w(x_1,1+t)= \lambda_3(\bar{\r},\bar{u}_1,\bar{\t})(x_1,t),\\
\di
\Sigma_3^{(j)}(\bar{\r},\bar{u}_1,\bar{\t})(x_1,t)=\Sigma_3^{(j)}(\rho_\pm,u_{1\pm},\t_\pm),\quad j=1,2,\quad
\bar u_2=\bar u_3=0,
\end{array} \right.\label{au}
\end{eqnarray}
where $\bar w(x_1,t)$ is the solution of Burgers equation $(\ref
{dbur})$ defined in \eqref{b-s}.
Then the planar 3-rarefaction wave $(\bar{\r},\bar{u},\bar{\t})(x_1,t)$ satisfies the Euler system
\begin{equation}\label{rare-s}
\left\{
\begin{array}{ll}
\di \bar\r_t+(\bar \r \bar u_1)_{x_1}=0,\\
\di (\bar \r \bar u_1)_t+(\bar \r\bar u_1^2+\bar p)_{x_1}=0,\\
\di  (\bar \r \bar u_i)_t+(\bar \r\bar u_1\bar u_i)_{x_1}=0,\qquad i=2,3,\\
\di \f{R}{\gamma-1}[(\bar\r\bar\t)_t+(\bar\r \bar u_1\bar\t)_{x_1}]+\bar p\bar u_{1x_1}=0
\end{array}
\right.
\end{equation}
with the initial values 
$$
(\bar\r_0,\bar u_0,\bar \t_0)(x_1):=(\bar\r,\bar u,\bar \t)(x_1,0)
$$
which is defined by using the smooth rarefaction wave for the Burgers equation evaluated at time $t=1$, as is suggested in \eqref{au}.
\bigskip

Now we can state the main result in this paper as follows.

\begin{theorem}\label{thm}
Let $(\bar\r, \bar u, \bar\t)(x_1,t)$ be the planar 3-rarefaction wave defined in \eqref{au}.
For each fixed state $(\r_+,u_+,\t_+)$, there exists a positive constant $\v_0$,
such that if $(\r_-,u_-,\t_-)\in R_3(\r_+,u_+,\t_+)$, and
\be
\v+\|(\r_0-\bar\r_0, u_0-\bar u_0,\t_0-\bar\t_0)\|_{H^2}\leq \v_0,
\ee
then the initial value problem \eqref{ns}--\eqref{far} admits a unique global smooth solution $(\r,u,\t)$ satisfying
\begin{equation}
\begin{cases}
(\r-\bar\r,u-\bar u,\t-\bar\t)\in C(0,+\infty;L^2(\Omega)),
\quad \na(\r,u,\t)\in C(0,+\infty;H^1(\Omega)),\\
\na^2\r\in L^2(0,+\infty;L^2(\Omega)),\quad \na^2(u,\t)\in L^2(0,+\infty; H^1(\Omega)),
\end{cases}
\end{equation}
and the time-asymptotic stability toward the planar rarefaction wave $(\bar\r, \bar u, \bar\t)(x_1,t)$ holds true:
\begin{equation}\label{large-time}
\lim_{t\rightarrow\infty} \sup_{x\in\Omega}|(\r,u,\t)(x,t)-(\bar\r,\bar u,\bar\t)(x_1,t)|=0.
\end{equation}

\begin{remark}
This is the first result about nonlinear stability of planar rarefaction wave for the three-dimensional non-isentropic equations,
while the corresponding stability results for shock wave or contact discontinuity are still completely open as far as we know.
\end{remark}

\begin{remark}
If we assume both $\|(\r_0-\r_0^r,u_0-u_0^r,\t_0-\t_0^r)\|_{L^2(\Omega)}+\|\na(\r_0,u_0,\t_0)\|_{H^1(\Omega)}$ and the wave strength $|(\r_+-\r_-,u_+-u_-,\t_+-\t_-)|$ are suitably small,
then the time-asymptotic stability of the 3-rarefaction wave fan  holds true:
$$
\lim_{t\rightarrow+\infty}\sup_{x\in\Omega}|(\r,u,\t)(x,t)-(\r^r,u^r,\t^r)(\f{x_1}{t})|=0,
$$
where $u_r=(u_1^r,0,0)^t$ and $(\r^r,u_1^r,\t^r)$ is the $3-$rarefaction wave to the Euler system \eqref{euler}, \eqref{R-in}.
\end{remark}
\end{theorem}

\bigskip

The rest part of the paper is arranged as follows. First, we present some properties on the smooth rarefaction wave solution in section \ref{Preliminaries}. 
Then, the energy estimates will be given in section \ref{Energy Estimates}. Finally, in the last section, based on
a priori estimates, we prove our main Theorem \ref{thm}.

%
%
%
%

\section{Rarefaction wave}\label{Preliminaries}
\setcounter{equation}{0}

In this section, we present some properties on the planar rarefaction wave constructed in \eqref{au}.

\begin{lemma}[\cite{MN-92}]\label{appr}
The problem~$(\ref{dbur})$ has a unique smooth global solution $\bar w(x_1,t)$ such that
\begin{itemize}
\item[(i)] $w_-<\bar w(x_1,t)<w_+, \  \bar w_{x_1}(x_1,t)>0,$
 \ for  $x_1\in\mathbb{R}, \ t\geq 0.$
\item[(ii)] For any $\ t> 0$ and p $\in[1,\infty]$, there exists a constant $C_{p,q}$ such that
\begin{equation*}
\begin{array}{ll}
\di \|\bar w(\cdot,t)-w^r(\frac\cdot t)\|_{L^p}\leq  C_{p,q}\v^{-\frac1p}(w_+-w_-), \\[3mm]
\di
\| \bar w_{x_1}(\cdot,t)\|_{L^p}\leq  C_{p,q}\min\{\v^{1-\frac{1}{p}}(w_+-w_-),
(w_+-w_-)^{\frac1p}t^{-1+\frac1p}\}, \\[3mm]
\di
  \| \bar w_{x_1x_1}(\cdot,t)\|_{L^p}\leq
C_{p,q}\min\{\v^{2-\frac1p}(w_+-w_-),\v^{(1-\frac{1}{2q})(1-\frac1p)}(w_+-w_-)^{-\frac{p-1}{2pq}}t^{-1-\frac{p-1}{2pq}}\},\\[3mm]
\di
\| \bar w_{x_1x_1x_1}(\cdot,t)\|_{L^p}\leq
C_{p,q}\min\{\v^{3-\frac1p}(w_+-w_-),\v^{(1-\f{1}{2q})(2-\f{1}{p})}(w_+-w_-)^{-\frac{2p-1}{2pq}}t^{-1-\frac{2p-1}{2pq}}\},\\[3mm]
\di
|\bar w_{x_1x_1}(x_1,t)|\leq C_{q}\v\bar w_{x_1}(x_1,t).
\end{array}
\end{equation*}
\item[(iii)] The smooth rarefaction wave $\bar w(x_1,t)$ and the original rarefaction wave $w^r(\frac{x_1}{t})$
are time-asymptotically equivalent, i.e.,
$$
\lim_{t\rightarrow+\infty}\sup_{x_1\in \mathbb{R}} |\bar w(x_1,t)-w^r(\f{x_1}{t})|=0.
$$
\end{itemize}
\end{lemma}

\begin{lemma}[\cite{JWX,MN-92}]\label{appu}
Let $\d=|(\r_+-\r_-, u_+-u_-, \t_+-\t_-)|$ is the strength of the smooth 3-rarefaction wave
$(\bar\r,\bar u, \bar\t)$ defined in (\ref{au}), then it satisfies the following properties:
\begin{itemize}
\item[(i)] $\bar u_{1x_1}(x_1,t)=\frac{2}{\gamma+1}\bar{w}_{x_1}>0,$
 \ for  $x_1\in\mathbb{R}, \ t\geq 0, \ \bar{\r}_{x_1} = \frac{1}{\sqrt{A\gamma \exp(\frac{\gamma - 1}{R}S_+)}}\bar{\r}^{\frac{3-\gamma}{2}}\bar{u}_{1x_1},\ \bar{\t}_{x_1} = \frac{\gamma - 1}{\sqrt{R\gamma}}\bar{\t}^{\f12}\bar{u}_{1x_1}.$
\item[(ii)] The following estimates hold for all $t> 0$ and
p $\in[1,\infty]$:
\begin{equation*}\begin{array}{l}
\di \|(\bar \r, \bar u_1,\bar \t)(\cdot,t)-(\r^r,u^r,\t^r)(\frac \cdot t)\|_{L^q}\leq C_{p,q}\delta\v^{-\frac1p}, \\
\di \| (\bar\r,\bar u_1,\bar\t)_{x_1}(\cdot,t)\|_{L^p}\leq
C_{p,q} \min\{\d\v^{1-\frac1p},\d^{\frac1p}(1+t)^{-1+\frac1p}\},\\
  \|(\bar\r,\bar u_1,\bar\t)_{x_1x_1}(\cdot,t)\|_{L^p}\leq
C_{p,q} \min\{\d\v^{2-\frac1p},\d^{-\frac{p-1}{2pq}}\v^{(1-\frac{1}{2q})(1-\frac1p)}\\
\qquad\qquad\qquad\qquad\qquad \times(1+t)^{-1-\frac{p-1}{2pq}}+\d^{\frac1p}(1+t)^{-2+\frac1p}\},\\
  \|(\bar\r,\bar u_1,\bar\t)_{x_1x_1x_1}(\cdot,t)\|_{L^p}\leq
C_{p,q} \min\{\d\v^{3-\frac1p},\d^{-\frac{2p-1}{2pq}}\v^{(1-\frac{1}{2q})(2-\frac1p)}\\
\qquad\qquad\qquad\qquad\qquad \times(1+t)^{-1-\frac{2p-1}{2pq}}+\d^{\frac1p}(1+t)^{-2+\frac1p}\}.
\end{array}\end{equation*}
\item[(iii)] Time-asymptotically, the smooth rarefaction wave and the inviscid rarefaction wave fan are equivalent, i.e.,
\begin{equation*}
\lim_{t\rightarrow+\infty}\sup_{x_1\in\mathbb{R}}| (\bar\r, \bar
u_1, \bar\t)(x_1,t)-(\r^r,u_1^r,\t^r)(\f{x_1}{t})|=0.
\end{equation*}
\end{itemize}
\end{lemma}

\

\

{\textbf{Notation.}} Throughout this paper, several positive generic constants are denoted by $C$ if without confusions.
For functional spaces, $H^s(\mathbb{R}\times \mathbb{T}^2)$ denotes the $s-$th order Sobolev space with the norm
$$
\|f\|_{H^s(\mathbb{R}\times \mathbb{T}^2)}\triangleq\sum_{j=0}^s\|\na_x^j f\| \quad {\rm and} \quad \|\cdot\|\triangleq\|\cdot\|_{L^2(\mathbb{R}\times \mathbb{T}^2)}.
$$

%
%
%
%

\section{A Priori Estimates}\label{Energy Estimates}
\setcounter{equation}{0}

\bigskip
Before we present the energy estimates, we first set
\be\label{perturb}
(\phi,\psi,\z)(x,t)=(\r-\bar\r, u-\bar u,\t-\bar\t)(x,t).
\ee
Then the solution is sought in the set of functional space $X(0,+\infty)$ defined by
$$
\ba
X(0,T)=\Big\{&(\phi,\psi,\z)|(\phi,\psi,\z)\in C(0,T;H^2),\quad\na \phi\in L^2(0,T;H^1),\\
& \na(\psi,\z)\in L^2(0,T;H^2)~{\rm and} ~\sup_{0\leq t\leq T}\|(\phi,\psi,\z)(t)\|_{H^2}\leq\chi \Big\},
\ea
$$
with $0\leq T\leq +\infty$.

Note that if $\chi$ is suitably small, then the condition $\di\sup_{0\leq t\leq T}\|(\phi,\psi,\z)\|_{H^2}\leq\chi$ and Sobolev embedding theorem imply that $|(\phi,\psi)|\leq \f12\r_-$, $|\z|\leq \f12\t_-$ and $|u|=|(u_1,u_2,u_3)|\leq C$ with $C$ being a positive constant which only depends on $\r_-,~u_{\pm}$. Therefore, the density function $\r(x,t):=\bar \r(x_1,t)+\phi(x,t)$ and the absolute temperature function $\t(x,t):=\bar\t(x_1,t)+\z(x,t)$ satisfy that
\be\label{upper-lower}
0<\f12\r_-\leq \r(x,t)\leq \f12\r_-+\r_+, \quad
0<\f12\t_-\leq \t(x,t)\leq \f12\t_-+\t_+,
\ee
since $0<\r_-\leq \bar\r(x_1,t)\leq\r_+$ and $0<\t_-\leq \bar\t(x_1,t)\leq\t_+$. It should be noted that the uniform lower and upper bounds of the density function $\r(x,t)$ in \eqref{upper-lower} guarantee the strict parabolicity of the momentum equation \eqref{ns}$_2$, which are crucial for the local and global-in-time existence of the classical solution to the system \eqref{ns}. Hence, for classical solutions, \eqref{ns} can be rewritten as
\begin{equation}\label{ns-1}
\begin{cases}
\di \r_t+u\cdot\na\r+\r\div u=0,\\
\di u_t+ u\cdot\na u +R\f{\t}{\r}\na\r+R\na\t =\f{1}{\r}(\mu\Delta u+(\mu+\lambda)\na \div u),\\
\di \f{R}{\gamma-1}(\t_t+  u\cdot\na\t)+R\t\div u =\f{1}{\r}\Big[\kappa\Delta\t+\frac{\mu}{2}|\na u + (\na u)^t|^2+\lambda(\div u)^2\Big],
\end{cases}
\end{equation}
with the initial data \eqref{initial} and far fields conditions on the $x_1$-direction \eqref{far}. From \eqref{rare-s} and \eqref{ns-1}, we can get the perturbation system for $(\phi,\psi,\z)$:
\begin{equation}\label{per}
\begin{cases}
\di \phi_t+u\cdot\na\phi+\r\div\psi+\psi\cdot\na\bar\r+\phi\div\bar u=0,\\
\di \psi_t+u\cdot\na\psi+R\f{\t}{\r}\na\phi+R\na\z+\psi\cdot\na\bar u+R\Big(\f{\t}{\r}-\f{\bar\t}{\bar\r}\Big)\na\bar\r\\
\di \quad =\f{1}{\r}\big(\mu\Delta\psi+(\mu+\lambda)\na\div \psi\big)+(\f{2\mu+\lambda}{\r}\bar u_{1x_1x_1},0,0)^t,\\
\di \f{R}{\gamma-1}(\z_t+u\cdot\na\z)+R\t\div\psi+\psi\cdot\na\bar\t+R\z\div \bar u=\f{\k}{\r}\Delta\z+\f{\k}{\r}\bar\t_{x_1x_1}\\
\di\quad +\f{1}{\r}\big[\frac{\mu}{2}|\na\psi + (\na\psi)^t|^2+ \lambda(\div \psi)^2+2\bar u_{1x_1}\big(2\mu\pa_1\psi_1+\lambda\div\psi\big)
+(2\mu+\lambda)\bar u_{1x_1}^2\big],
\end{cases}
\end{equation}
and the initial data is
\be\label{per-data}
(\phi,\psi,\z)(x,0)=(\phi_0,\psi_0,\z_0)(x)=(\r_0-\bar\r_0,u_0-\bar u_0,\t_0-\bar\t_0)(x).
\ee

\

Since the proof for the local-in-time existence and uniqueness of the classical solution to \eqref{per}-\eqref{per-data} is standard (for instance, one can refer to \cite{Nash} or \cite{Solo}), in particular for the suitably small perturbation of the solution around the planar rarefaction wave satisfying the property \eqref{upper-lower}, the details will be omitted. To prove Theorem \ref{thm}, it suffices to show the following a priori estimates.

\begin{proposition}\label{priori}
{\rm (A priori estimates)} Suppose that the reformulated problem \eqref{per}-\eqref{per-data} admits a solution $(\phi,\psi,\z)\in X(0,T)$ for some $T>0$. Then there exist positive constants $\chi\leqq1$ and $C$ independent of $T$, such that if
\begin{equation}\label{assumption}
\ba
\sup_{0\leq t \leq
T}\|(\phi,\psi,\z)(\cdot,t)\|_{H^2}\leq \chi,
\ea
\end{equation}
then it follows the estimates:
\begin{equation}\label{full-es}
\begin{array}{l}
\di \sup_{0\leq t\leq T}\|(\phi,\psi,\z)(\cdot,t)\|_{H^2}^2+\int_0^T\Big[\|\sqrt{\bar u_{1x_1}}(\phi,\psi,\zeta)\|^2+\|\na\phi\|^2_{H^1}+\|\na(\psi,\z)\|^2_{H^2}\Big]d\tau\\
\di\leq C(\|(\phi_0,\psi_0,\z_0)\|^2_{H^2}+\v^{\frac18}).
\end{array}
\end{equation}

\end{proposition}

\

From now on, we always assume that $\chi+\v\leqq 1$. Proposition \ref{priori} is an easy consequence of the following lemmas. We first give the following $L^2$ estimate.


\begin{lemma} \label{le-1}
For $T>0$ and $(\phi,\psi,\z)\in X(0,T)$ satisfying a priori assumption \eqref{assumption} with suitably small $\chi+\v$, we have for $t\in[0,T]$,
\begin{equation}\label{basic}
\begin{array}{ll}
\di \|(\phi,\psi,\z)(t)\|^2+\int_0^t\big[\|\sqrt{\bar u_{1x_1}}(\phi,\psi_1,\z)\|^2+\|\na(\psi,\z)\|^2\big]d\tau
 \leq C\|(\phi_0,\psi_0,\z_0)\|^2+C\v^{\f18}.
\end{array}
\end{equation}

\end{lemma}

{\textbf{Proof}:}
For ideal polytropic fluids, it holds
$$
 S=-R\ln\rho+ \frac{R}{\gamma-1}\ln\t+\f{R}{\gamma-1}\ln\f{R}{A},\\
\qquad  p=R\rho\theta=A\rho^{\gamma}\exp\Big(\f{\gamma-1}{R}S\Big).
$$
Denote
$$
\begin{array}{l}
\displaystyle \mathbf{X}=\Big(\rho,\r u_1,\r u_2,\r u_3,\r \Big(\f{R}{\gamma-1}\t+\f{|u|^2}{2}\Big)\Big)^t,\\
\displaystyle \mathbf{Y}=\Big(\r u,\r u u_1+p\mathbb{I}_1,\r uu_2+p\mathbb{I}_2,\r u u_3+p\mathbb{I}_3,\r u\Big(\f{R}{\gamma-1}\t+\f{|u|^2}{2}\Big)+pu\Big)^t,
\end{array}
$$
where $\mathbb{I}_1=(1,0,0)^t, \mathbb{I}_2=(0,1,0)^t, \mathbb{I}_3=(0,0,1)^t.$
Then the system \eqref{ns} can be rewritten as
$$
\mathbf{X}_t+\div\mathbf{Y}=
\left(
\begin{array}{c}
0\\
\displaystyle \mu\Delta u_1+(\mu+\lambda)\pa_1\div u\\
\displaystyle \mu\Delta u_2+(\mu+\lambda)\pa_2\div u\\
\displaystyle \mu\Delta u_3+(\mu+\lambda)\pa_3\div u\\
\displaystyle \k\Delta\t+\div(u\mathcal{T})
\end{array}
\right),
$$
where $\pa_j=\pa_{x_j}~(j=1,2,3)$. We define a relative entropy-entropy flux pair $(\eta,q)$ as
$$
\left\{
\begin{array}{l}
\displaystyle \eta=\bar\theta\left\{-\rho
S+\bar\rho\bar{S}+\nabla_\mathbf{X}(\rho
S)\Big|_{\mathbf{X}=\bar{\mathbf{X}}}\cdot(\mathbf{X}-\bar{\mathbf{X}})\right\},\\
\displaystyle q_j=\bar\theta\left\{-\rho u_j
S+\bar\rho\bar{u}_j\bar{S}+\nabla_\mathbf{X}(\rho
S)\Big|_{\mathbf{X}=\bar{\mathbf{X}}}\cdot(\mathbf{Y}_j-\bar{\mathbf{Y}}_j)\right\}\quad j=1,2,3.
\end{array}
\right.
$$
Here, we can compute that
$$
\displaystyle (\rho S)_{\rho}=S+\frac{|u|^2}{2\theta}-\frac{R\gamma}{\gamma-1},\qquad
\displaystyle (\rho S)_{m_i}=-\frac{u_i}{\theta},~i=1,2,3,\qquad
\displaystyle (\rho S)_{\mathcal{E}}=\frac{1}{\theta},
$$
where $m_i=\r u_i~(i=1,2,3)$ and $\mathcal{E}=\r(\f{R}{\gamma-1}\t+\f{|u|^2}{2})$, then
$$
\left\{\begin{aligned}
 \eta&=\f{R}{\gamma-1}\rho\theta-\bar\theta\rho
S+\rho\Big[\Big(\bar{S}-\frac{R\gamma}{\gamma-1}\Big)\bar\theta+\frac{|u-\bar u|^2}{2}\Big]
+R\bar\rho\bar\theta \\
&=R\rho\bar\t\Psi\left(\f{\bar\r}{\rho}\right)+\f{R}{\gamma-1}\rho\bar\t\Psi\left(\f{\t}{\bar\t}\right)
+\f12\rho|u-\bar u|^2,\\
q&=u\eta +R(u-\bar u)(\rho\theta-\bar\rho\bar\theta),
\end{aligned}
\right.
$$
where $\Psi(\cdot)$ is the convex function
$$
\Psi(s)=s-\ln s-1.
$$
Then, for $\mathbf{X}$ in any closed bounded region in
$\sum=\{\mathbf{X}:\rho>0,\theta>0\}$, there exists a positive
constant $C_0$ such that
$$
C_0^{-1}|(\phi,\psi,\z)|^2\leq\eta\leq
C_0|(\phi,\psi,\z)|^2.
$$
Direct computations yield that
\begin{equation}\label{E1}
\begin{array}{l}
\di\eta_t+\div q+\f{\bar\t}{\t}\big(\frac{\mu}{2}|\na\psi + (\na\psi)^t|^2+ \lambda(\div \psi)^2\big)+\f{\k\bar\t}{\t^2}|\na\z|^2
-\big[\nabla_{(\bar\rho,\bar{u},\bar{S})}\eta\cdot(\bar\rho,\bar{u},\bar{S})_t
+\nabla_{(\bar\rho,\bar{u},\bar{S})}q\cdot(\bar\rho,\bar{u},\bar{S})_{x_1}\big]\\
\di=\div\Big[\psi\big(\mu\na u+(\mu+\lambda)\div u\big)+\f{\k\z\na\z}{\t}\Big]
-\pa_1(\mu\psi_1\bar u_{1x_1})-\div\big[(\mu+\lambda)\psi\bar u_{1x_1}\big]\\
\di - \mu|\na\psi|^2 - (\mu + \lambda)(\div\psi)^2 + \frac{\mu}{2}|\na\psi + (\na\psi)^t|^2 + \lambda(\div\psi)^2\\
\di +\f{2\z}{\t}\big(2\mu\pa_1\psi_1+\lambda\div\psi\big)\bar u_{1x_1}
+\f{\k}{\t^2}\z\pa_1\z\bar\t_{x_1}+(2\mu+\lambda)\big(\f{\z}{\t}\bar u^2_{1x_1}+\psi_1\bar u_{1x_1x_1}\big)+\f{\k}{\t}\z\bar\t_{x_1x_1}.
\end{array}
\end{equation}
There exists a positive constant $C>0$ such that (cf. \cite{LWW})
$$
\begin{array}{ll}
\di\quad -\big[\nabla_{(\bar\rho,\bar{u},\bar{S})}\eta\cdot(\bar\rho,\bar{u},\bar{S})_t
+\nabla_{(\bar\rho,\bar{u},\bar{S})}q\cdot(\bar\rho,\bar{u},\bar{S})_{x_1}\big]\\[2mm]
\di =\bar u_{1x_1}\Big[\rho\psi_1^2+R(\gamma-1)\rho\bar\t \Psi\Big(\f{\bar\rho}\rho\Big)
+R\rho\bar\t\Psi\Big(\f{\t}{\bar\t}\Big)\Big]
+\bar\t_{x_1}\rho\psi_1\Big(R\ln\f{\bar\r}{\r}+\f{R}{\gamma-1}\ln\f{\t}{\bar\t}\Big)\\
\di  \geq C^{-1} \bar u_{1x_1}(\phi^2+\psi_1^2+\z^2).
\end{array}
$$
Integrating \eqref{E1} with respect to $x,t$ over $\Omega\times(0,t)$ yields that
\begin{equation}\label{E2}
\ba
&\quad \|(\phi,\psi,\z)(t)\|^2+\int_0^t\Big[\|(\na\psi,\na\z)\|^2+
\|\sqrt{\bar u_{1x_1}}(\phi,\psi_1,\z)\|^2\Big]d\tau
\\
& \leq C\|(\phi_0,\psi_0,\z_0)\|^2+C\Big|\int_0^t\int
\Big[\f{2\z}{\t}\big(2\mu\pa_1\psi_1+\lambda\div\psi\big)\bar u_{1x_1}
+\f{\k}{\t^2}\z\pa_1\z\bar\t_{x_1}\Big]dxd\tau\Big|\\
&+C\Big|\int_0^t\int\Big[\psi_1\bar u_{1x_1x_1}+\f{\z}{\t}\bar u^2_{1x_1}+\f{\k}{\t}\z\bar\t_{x_1x_1}\Big] dxd\tau\Big|,
\ea
\end{equation}
where we have used the following fact
\[\int \mu|\na\psi|^2 + (\mu + \lambda)(\div\psi)^2dx = \int \frac{\mu}{2}|\na\psi + (\na\psi)^t|^2 + \lambda(\div\psi)^2dx.\]
First, by the Cauchy's inequality and Lemma \ref{appu}, it holds that
\be\label{e1}
\begin{array}{l}
\di C\Big|\int_0^t\int
\Big[\f{2\z}{\t}\big(2\mu\pa_1\psi_1+\lambda\div\psi\big)\bar u_{1x_1}
+\f{\k}{\t^2}\z\pa_1\z\bar\t_{x_1}\Big]dxd\tau\Big|\\
\di\leq \f12\int_0^t\|\na(\psi,\z)\|^2d\tau+C\v\int_0^t\|\sqrt{\bar u_{1x_1}}\z\|^2d\tau.
\end{array}
\ee
By Sobolev's inequality, H\"{o}lder's inequality, Young's inequality, Lemma \ref{appu} and assumption \eqref{assumption}, we have
\begin{equation}\label{e2}
\begin{array}{ll}
\di C\Big|\int_0^t\int \psi_1 \bar u_{1x_1x_1}dxd\tau\Big|
\leq C\int_0^t\int_{\mathbb{T}^2}\|\psi_1\|_{L^{\infty}(\mathbb{R})}\|\bar u_{1x_1x_1}\|_{L^1(\mathbb{R})} dx_2dx_3d\tau\\
\di \leq C\v^{\f18}\int_0^t(1+\tau)^{-\f78}\Big(\int_{\mathbb{T}^2}\|\psi_1\|^{\f12}_{L^2(\mathbb{R})}
\|\pa_1\psi_1\|^{\f12}_{L^2(\mathbb{R})}dx_2dx_3\Big)d\tau\\
\di \leq C\v^{\f18}\int_0^t(1+\tau)^{-\f78}\|\pa_1\psi_1\|^{\f12}\Big(\int_{\mathbb{T}^2}\|\psi_1\|^{\f23}_{L^2(\mathbb{R})}dx_2dx_3\Big)^{\f34}d\tau\\
\di \leq C\v^{\f18}\int_0^t(1+\tau)^{-\f78}\|\psi_1\|^{\f12}\|\pa_1\psi_1\|^{\f12}d\tau
 \leq C\v^{\f18}\int_0^t\|\pa_1\psi_1\|^2d\tau+C\v^{\f18}.
\end{array}
\end{equation}
Similarly, one has
\be
\begin{array}{l}\label{e3}
\di C\Big|\int_0^t\int\Big[\f{\z}{\t}\bar u^2_{1x_1}+\f{\k}{\t}\z\bar\t_{x_1x_1}\Big] dxd\tau\Big|
\leq C\v^{\f18}\int_0^t\|\na\z\|^2d\tau+C\v^{\f18}.
\end{array}
\ee

Substituting the estimates \eqref{e1}-\eqref{e3} into \eqref{E2} gives \eqref{basic}, and the proof of Lemma \ref{le-1} is completed.

\hfill $\Box$

Next, we want to get the estimation of $\na\phi$. Compared with the one-dimensional stability results in \cite{MN-86,MN-92,NYZ}, the physical viscosity
in momentum equation \eqref{per}$_2$ has the form: $\mu\Delta \psi+(\mu+\lambda)\na \div \psi$ in high dimensions. Therefore, we can not substitute mass equation \eqref{per}$_1$ into momentum equation \eqref{per}$_2$ to obtain the derivative estimate of density perturbation $\na\phi$ as in \cite{MN-86,MN-92,NYZ}. Our new observation is that we find some cancellation between the flux terms and viscosity terms for system \eqref{per}, by which we successfully overcome the difficulty when the planar rarefaction wave propagate in $x_2,x_3$-directions may interact with $x_1$-direction, and derive the derivative estimates of density perturbation $\na\phi$. The following lemma is crucial to get a priori estimates \eqref{full-es}.
\begin{lemma}\label{le-2}
For $T>0$ and $(\phi,\psi,\z)\in X(0,T)$ satisfying a priori assumption \eqref{assumption} with suitably small $\chi+\v$, it holds that for $t\in[0,T]$,
\be\label{phi-x}
\di \|\na\phi(t)\|^2+\int_0^t\|\na\phi\|^2d\tau\leq C\|(\phi_0,\psi_0,\z_0,\na\phi_0)\|^2+C\v^{\f18}+C(\chi+\v)\int_0^t\|\na^3 u\|^2d\tau.
\ee
\end{lemma}

{\textbf{Proof}:} For this, we multiply \eqref{per}$_2$ by $\r\na\phi$ and integrate by parts with respect to $x$ to obtain
\begin{equation}\label{phi-11}
\begin{array}{ll}
\di \int\r\psi_t\cdot\na\phi dx+\int\r u\cdot \nabla \psi\cdot \nabla \phi dx+\int R\t|\na\phi|^2dxd\tau\\[2mm]
\di =-\int\rho\Big[R\na\z+\psi\cdot\na\bar u
+R\Big(\f{\t}{\r}-\f{\bar\t}{\bar\r}\Big)\na\bar\r\Big]\cdot\na\phi\,dx\\
\di\quad+\int \Big[\mu\Delta \psi+(\mu+\lambda)\na\div \psi\Big]\cdot\nabla\phi\,dx+(2\mu+\lambda)\int \bar u_{1x_1x_1}\pa_{1}\phi \,dx
,
\end{array}
\end{equation}
By using the following two facts:
\be
\begin{array}{ll}
\di \int\r\psi_t\cdot\na\phi dx+\int\r u\cdot \nabla \psi\cdot \nabla \phi dx\\
\di=\f{d}{dt}\int \r\psi\cdot\na\phi \,dx+\int\div(\r u)\psi\cdot\na\phi\,dx
+\int \div(\r\psi)\phi_t \,dx+\int\r u\cdot \nabla \psi\cdot \nabla \phi dx\\
\di=\f{d}{dt}\int \r\psi\cdot\na\phi \,dx+\int\partial_{j}(\r u_j\psi)\cdot\na\phi\,dx -\int\div(\r\psi)(u\cdot\na\phi+\r\div\psi+\psi\cdot\na\bar\r+\phi\div\bar u)dx
\end{array}
\ee
and
\be
\di \int\Big[\mu\Delta\psi+(\mu+\lambda)\na\div\psi\Big]\cdot\na\phi\,dx
=(2\mu+\lambda)\int\na\phi\cdot\na\div\psi \, dx,
\ee
the equality \eqref{phi-11} becomes
\begin{equation}\label{phi-1}
\begin{array}{ll}
\di \f{d}{dt}\int\r\psi\cdot\na\phi dx+\int R\t|\na\phi|^2dxd\tau
=\int\div(\r\psi)(u\cdot\na\phi+\r\div\psi+\psi\cdot\na\bar\r+\phi\div\bar u)dx\\[2mm]
\di\quad -\int\partial_{j}(\r u_j\psi)\cdot\na\phi\,dx-\int\r\Big[R\na\z+\psi\cdot\na\bar u
+R\Big(\f{\t}{\r}-\f{\bar\t}{\bar\r}\Big)\na\bar\r\Big]\cdot\na\phi\,dx\\
\di\quad
+(2\mu+\lambda)\int \na\phi\cdot\na\div \psi\,dx+(2\mu+\lambda)\int \bar u_{1x_1x_1}\pa_1\phi \,dx.
\end{array}
\end{equation}
In order to close the a priori assumption \eqref{assumption}, we need to get rid of the higher order term $\di (2\mu+\lambda)\int \na\phi\cdot\na\div \psi\,dx$ in \eqref{phi-1}. Otherwise, the first-order derivative estimate  in \eqref{phi-1} will depend on the second order derivative $\na\div \psi$ and deductively one can not close the a priori assumption \eqref{assumption}. For this, we first apply $\pa_i~(i=1,2,3)$ to the equation \eqref{per}$_1$ to derive
\be\label{phi-2}
\begin{array}{l}
\di \pa_i\phi_t+u\cdot\na\pa_i\phi+\r\pa_i\div\psi+\pa_i u\cdot\na\phi+\pa_i\r\div\psi+\pa_i\psi\cdot\na\bar\r
+\pa_i\phi\div \bar u+\psi\cdot\na\pa_i\bar\r+
\phi\pa_i\div\bar u=0.
\end{array}
\ee
Then multiplying the above equation by $\f{2\mu+\lambda}{\r}\pa_i\phi$,  integrating over the domain $\Omega$ with respect to $x$ 
and summing $i$ from $1$ to $3$ yield
\be\label{phi-3}
\begin{array}{l}
\di \f{d}{dt}\int\f{2\mu+\lambda}{2\r}|\na\phi|^2dx
=(2\mu+\lambda)\int\f{\div u}{\r}|\na\phi|^2dx-\int\f{2\mu+\lambda}{\r}\pa_i\phi\big(
\pa_i u\cdot\na\phi+\pa_i\r\div\psi\\
\di\quad +\pa_i\psi\cdot\na\bar\r+\pa_i\phi\div \bar u+\psi\cdot\na\pa_i\bar\r+\phi\pa_i\div\bar u\big)dx-(2\mu+\lambda)\int\na\phi\cdot\na\div\psi \,dx,
\end{array}
\ee
where we have the following equality:
$$
\ba
\int\f{2\mu+\lambda}{\r}\pa_i\phi(\pa_i\phi_t+u\cdot\na\pa_i\phi)
&=\f{d}{dt}\int\f{2\mu+\lambda}{2\r}|\pa_i\phi|^2dx
-\f{2\mu+\lambda}{2}\int\Big[\Big(\f{1}{\r}\Big)_t+\div\Big(\f{u}{\r}\Big)\Big]|\pa_i\phi|^2dx\\
&=\f{d}{dt}\int\f{2\mu+\lambda}{2\r}|\pa_i\phi|^2dx-(2\mu+\lambda)\int\f{\div u}{\r}|\pa_i\phi|^2dx.
\ea
$$
Thus adding the equalities \eqref{phi-1} and  \eqref{phi-3} together and the higher order term $\di (2\mu+\lambda)\int \na\phi\cdot\na\div \psi\,dx$ will be cancelled as desired, and then integrating the resulted equation with respect to the time $t$ over $(0,t)$ to give
\be\label{phi-4}
\begin{array}{l}
\di \int\Big(\f{2\mu+\lambda}{2\r}|\na\phi|^2+\r\psi\cdot\na\phi\Big)dx\Big|_{\tau=0}^{\tau=t}
+\int_0^t\int R\t|\na\phi|^2dxd\tau\\
\di=(2\mu+\lambda)\int_0^t\int\f{\div u}{\r}|\na\phi|^2dxd\tau-\int_0^t\int\f{2\mu+\lambda}{\r}\pa_i\phi\big(
\pa_i u\cdot\na\phi+\pa_i\r\div\psi\\
\di\quad+\pa_i\psi\cdot\na\bar\r+\pa_i\phi\div \bar u+\psi\cdot\na\pa_i\bar\r+\phi\pa_i\div\bar u\big)dxd\tau
-\int_0^t\int\partial_{j}(\r u_j\psi)\cdot\na\phi\,dxd\tau\\
\di\quad +\int_0^t\int\div(\r\psi)(u\cdot\na\phi+\r\div\psi+\psi\cdot\na\bar\r+\phi\div\bar u)dxd\tau
 \\
\di\quad-\int_0^t\int\r\Big[R\na\z+\psi\cdot\na\bar u
+R\Big(\f{\t}{\r}-\f{\bar\t}{\bar\r}\Big)\na\bar\r\Big]\cdot\na\phi\,dxd\tau
+(2\mu+\lambda)\int_0^t\int \bar u_{1x_1x_1}\pa_1\phi \,dxd\tau.
\end{array}
\ee
We just estimate the first term on the right hand side of \eqref{phi-4} as follows and the other terms can be done similarly and the details for estimating these terms will be omitted for brevity. By Sobolev's inequality Lemma \ref{appu} and assumption \eqref{assumption}, one has
\be\label{hard}
\begin{array}{l}
\di(2\mu+\lambda)\int_0^t\int\f{\div u}{\r}|\na\phi|^2dxd\tau\leq C\int_0^t\int|\div u||\na\phi|^2 dxd\tau\\
\di \leq \int_0^t\|\na u\|_{L^{\infty}}\|\na\phi\|^2d\tau\leq C\int_0^t\|\na u\|_{H^2}\|\na\phi\|^2d\tau\\
\di\leq C\int_0^t(\|\na u\|+\|\na^2 u\|)\|\na\phi\|^2d\tau+C\int_0^t\|\na^3 u\|\|\na\phi\|^2d\tau\\
\di \leq C\sup_{0\leq\tau\leq t}\|(\na u,\na^2 u)(\tau)\|\int_0^t\|\na\phi\|^2d\tau
+C\sup_{0\leq\tau\leq t}\|\na \phi(\tau)\|\int_0^t\|\na^3u\|\|\na\phi\|d\tau\\
\di\leq C(\chi+\v)\int_0^t(\|\na\phi\|^2+\|\na^3 u\|^2)d\tau.
\end{array}
\ee
By Cauchy's inequality and the estimates as in \eqref{hard}, it follows from \eqref{phi-4} that
\begin{equation}
\begin{array}{ll}
\di \|\na\phi(t)\|^2 +\int_0^t\|\na\phi\|^2d\tau\leq C\|(\psi_0,\na\phi_0)\|^2+C\v+C\|\psi(t)\|^2+C\int_0^t\|\na(\psi,\z)\|^2d\tau\\
\di\quad\quad\quad\quad\quad\quad\quad\quad\quad\quad\quad
 +C\v\int_0^t\|\sqrt{\bar u_{1x_1}}(\phi,\psi_1,\z)\|^2d\tau+C(\chi+\v)\int_0^t\|\na^3 u\|^2d\tau,
\end{array}
\end{equation}
which together with \eqref{basic} leads to \eqref{phi-x}, and the proof of Lemma \ref{le-2} is completed.

\hfill $\Box$

\

\begin{lemma}\label{high-1}
For $T>0$ and $(\phi,\psi,\z)\in X(0,T)$ satisfying a priori assumption \eqref{assumption} with suitably small $\chi+\v$, we have for $t\in[0,T]$,
\be\label{psi-1x}
\begin{array}{l}
\di \|\na(\psi,\z)(t)\|^2+\int_0^t\|\na^2(\psi,\z)\|^2d\tau
\leq C\|(\phi_0,\psi_0,\z_0)\|^2_{H^1}+C\v^{\f18}+C(\chi+\v)\int_0^t\|\na^3 u\|^2d\tau.
\end{array}\ee
\end{lemma}

\

{\textbf {Proof}:}  Multiplying the equation \eqref{per}$_2$ by $(-\Delta\psi)$, and integrating over $\Omega\times(0,t)$ lead to
\be\label{psi-2}
\begin{array}{l}
\di\quad \int\f{|\na\psi|^2}{2}dx\Big|_{\tau=0}^{\tau=t}+\int_0^t\int\Big(\f{\mu}{\r}|\Delta\psi|^2+\f{\mu+\lambda}{\r}|\na\div \psi|^2\Big)dxd\tau\\
\di=\int_0^t\int \Big[u\cdot\na\psi+R\f{\t}{\r}\na\phi+R\na\z+\psi\cdot\na\bar u+R\Big(\f{\t}{\r}-\f{\bar\t}{\bar\r}\Big)\na\bar\r\Big]\Delta\psi \,dxd\tau\\
\di +\int_0^t\int\f{(\mu+\lambda)}{\r^2}\big(
\pa_j\r\pa_j\div\psi\div\psi -\pa_i\r\div\psi\Delta\psi_i\big)dxd\tau\\
\di
-(2\mu+\lambda)\int_0^t\int\f{1}{\r}\bar u_{1x_1x_1}\Delta\psi_1\,dxd\tau:=\sum_{i=1}^3I_i.
\end{array}
\ee
Here we use the following fact:
\be\label{psi-3}
\begin{array}{l}
\di \int\f{1}{\r}\big(\mu+\lambda\big)\na\div\psi\Delta\psi \,dx
=\int\f{\mu+\lambda}{\r}|\na\div \psi|^2dx
-\int\f{\mu+\lambda}{\r^2}(\pa_j\r\pa_j\div\psi\div\psi-\pa_i\r\div\psi\Delta\psi_i)\, dx.
\end{array}
\ee

Now we will estimate each $I_i~(i=1,2,3)$ on the right hand side of \eqref{psi-2}. By Cauchy's inequality and Lemma \ref{appu}, one has
\be\label{I1}
\di |I_1|\leq \s\int_0^t\|\na^2\psi\|^2d\tau+C_{\s}\int_0^t\|\na(\phi,\psi,\z)\|^2d\tau
+C_{\s}\v\int_0^t\|\sqrt{\bar u_{1x_1}}(\phi,\psi_1,\z)\|^2d\tau,
\ee
where $\s$ is a suitably small positive constant to be determined and $C_\s$ is a positive constant depending on $\s$. It follows from Cauchy's inequality,
Sobolev's inequality, Lemma \ref{appu} and assumption \eqref{assumption} that
\be\label{I2}
\ba
|I_2|&\leq C\int_0^t\int (|\na\phi|+|\bar\r_{x_1}|)|\na\psi||\na^2\psi|dxd\tau\\
& \leq C \int_0^t\|\na\phi\|_{L^4}\|\na\psi\|_{L^4}\|\na^2\psi\|d\tau
+C\v\int_0^t\|\na\psi\|\|\na^2\psi\|d\tau\\
& \leq C\int_0^t\|\na\phi\|_{H^1}\|\na\psi\|_{H^1}\|\na^2\psi\|d\tau
+C\v\int_0^t\|\na\psi\|\|\na^2\psi\|d\tau\\
&\leq C(\chi+\v)\int_0^t\|(\na\psi,\na^2\psi)\|^2d\tau.
\ea
\ee
By Cauchy's inequality and Lemma \ref{appu}, we obtain
\be\label{I3}
\di |I_3|\leq \s\int_0^t\|\na^2\psi\|^2d\tau+C_{\s}\int_0^t\|\bar u_{1x_1x_1}\|^2d\tau
\leq \s\int_0^t\|\na^2\psi\|^2d\tau+C_{\s}\v.
\ee
Substituting \eqref{I1}-\eqref{I3} into \eqref{psi-2} yields
\be\label{psi-4}
\begin{array}{l}
\di\quad \|\na\psi(t)\|^2+\int_0^t\|\na^2\psi\|^2d\tau\\
\di\leq C\|\na\psi_0\|^2+C\v+C\int_0^t\|\na(\phi,\psi,\z)\|^2d\tau
+C\v\int_0^t\|\sqrt{\bar u_{1x_1}}(\phi,\psi_1,\z)\|^2d\tau.
\end{array}
\ee

Next, we  estimate $\|\na\z\|$. We multiply the equation \eqref{per}$_3$ by $(-\Delta\z)$,
and integrate by parts over $\Omega\times(0,t)$, similar as \eqref{psi-3}, it holds
\be\label{wt-1}
\begin{array}{l}
\di \f{R}{\gamma-1}\int\f{|\na\z|^2}{2}dx\Big|_{\tau=0}^{\tau=t}+\int_0^t\int\f{\k}{\r}|\na^2\z|^2dxd\tau
=\int_0^t\int\Big[\f{R}{\gamma-1}u\cdot\na\z+R\t\div\psi\\
\di\quad+\psi\cdot\na\bar\t+R\z\div\bar u\Big]\Delta\z\,dxd\tau-\int_0^t\int\f{\Delta\z}{\r}\big(\frac{\mu}{2}|\na\psi + (\na\psi)^t|^2 + \lambda(\div\psi)^2\big)dxd\tau\\
\di\quad -\int_0^t\int\f{2\Delta\z}{\r}\bar u_{1x_1}\big(2\mu\pa_1\psi_1+\lambda\div\psi\big)dxd\tau
-\int_0^t\int\Delta\z\Big(\f{\k}{\r}\bar\t_{x_1x_1}+(2\mu+\lambda)\bar u_{1x_1}^2\Big)dxd\tau\\
\di\quad+\int_0^t\int\f{\k}{\r^2}(\pa_i\r\pa_i\pa_j\z\pa_j\z-\pa_j\r\pa_j\z\Delta\z)dxd\tau:=\sum_{i=4}^8I_i.
\end{array}
\ee
We just estimate $I_5$, other terms are similar to \eqref{I1}-\eqref{I3}. By H\"older's inequality, Sobolev's inequality, Cauchy's inequality, Lemma \ref{appu} and assumption \eqref{assumption},
\be\label{wt-2}
\ba
|I_5|&\leq C\int_0^t\|\na^2\z\|\|\na\psi\|^2_{L^4}d\tau\leq C\int_0^t\|\na^2\z\|\|\na\psi\|_{H^1}^2d\tau\\
&\leq \sigma\int_0^t\|\na^2\z\|^2d\tau+C_\sigma\int_0^t\|(\na\psi,\na^2\psi)\|^2d\tau.
\ea
\ee
Similar to \eqref{psi-4}, it follows from \eqref{wt-1} and \eqref{wt-2} that
\be\label{zeta-2}
\begin{array}{l}
\di\|\na\z(t)\|^2+\int_0^t\|\na^2\z\|^2d\tau \leq C\|\na\z_0\|^2+C\v+C\int_0^t\|\na(\psi,\z)\|^2d\tau\\
\di\quad\quad\quad\quad
+C\v\int_0^t\|\sqrt{\bar u_{1x_1}}(\psi_1,\z)\|^2d\tau+C\int_0^t\|\na^2\psi\|^2d\tau,
\end{array}
\ee
which together with \eqref{psi-4} leads to
\be\label{zeta-3}
\begin{array}{l}
\di\quad \|\na(\psi,\z)(t)\|^2+\int_0^t\|\na^2(\psi,\z)\|^2d\tau\\
\di\leq C\|\na(\psi_0,\z_0)\|^2+C\v+C\int_0^t\|\na(\phi,\psi,\z)\|^2d\tau
+C\v\int_0^t\|\sqrt{\bar u_{1x_1}}(\phi,\psi_1,\z)\|^2d\tau.
\end{array}
\ee
Finally, combining \eqref{zeta-3} with \eqref{basic} and \eqref{phi-x} implies \eqref{psi-1x}, and the proof of Lemma \ref{high-1} is completed.

\hfill $\Box$

\

The following lemmas are concerned with the higher order estimates of the perturbation $(\phi,\psi,\z)$. In order to obtain these estimates,
we prefer to consider the system \eqref{ns-1} of $(\r,u,\t)$ rather than the perturbation system \eqref{per} of $(\phi,\psi,\z)$ due to the fact
$\|(\bar\r,\bar u_1,\bar\t)_{x_1x_1}, (\bar\r,\bar u_1,\bar\t)_{x_1x_1x_1}\|^2\backsim \v(1+t)^{-2}$, which is integrable with respect to the time $t$ on $\mathbb{R}^+$. Therefore, we can use system \eqref{ns-1} to
derive Lemmas \ref{le-3} and \ref{le-5}. We start from Lemma \ref{le-3} concerning the second order derivatives estimates for $\phi$.

\begin{lemma}\label{le-3}
For $T>0$ and $(\phi,\psi,\z)\in X(0,T)$ satisfying a priori assumption \eqref{assumption} with suitably small $\chi+\v$, it holds that for $t\in[0,T]$,
\be\label{phi-xx}
\di \|\na^2\phi(t)\|^2+\int_0^t\|\na^2\phi\|^2d\tau\leq C(\|(\phi_0,\psi_0,\z_0)\|_{H^1}^2+\|\na^2\phi_0\|^2+\v^{\f18})
+C(\chi+\v)\int_0^t\|\na^3 u\|^2d\tau.
\ee
\end{lemma}

{\textbf{Proof}:} Applying $\pa_j\pa_i~(i,j=1,2,3)$ to the mass equation \eqref{ns-1}$_1$ and
$\pa_j~(j=1,2,3)$ to the $i-$th $(i=1,2,3)$ component of the momentum equation \eqref{ns-1}$_2$, we have
\be\label{rho-x1}
\begin{cases}
\di \pa_j\pa_i\r_t+u\cdot\na\pa_j\pa_i\r+\r\pa_j\pa_i\div u+(\pa_j u\cdot\na\pa_i\r+\pa_i u\cdot\na\pa_j\r+\div u\pa_j\pa_i\r)\\
\di\quad+(\pa_j\r\pa_i\div u+\pa_j\pa_i u\cdot\na\r+\pa_i\r\pa_j\div u)=0,\\
\di \pa_j u_{it}+u\cdot\na\pa_j u_i+R\f{\t}{\r}\pa_j\pa_i\r+R\pa_j\pa_i\t+\pa_j u\cdot\na u_i+\pa_j\Big(R\f{\t}{\r}\Big)\pa_i\r\\
\di\quad=\f{1}{\r}\big(\mu\Delta\pa_j u_i+(\mu+\lambda)\pa_j\pa_i\div u\big)-\f{\pa_j\r}{\r^2}\big(\mu\Delta u_i+(\mu+\lambda)\pa_i\div u\big).
\end{cases}
\ee

Next, multiplying the equation \eqref{rho-x1}$_2$ by $\r\pa_j\pa_i\r$ and integrating with respect to $x$ lead to
\begin{equation}\label{rho-x2}
\begin{array}{ll}
\di \f{d}{dt}\int\r\pa_j u_i\pa_j\pa_i\r dx+\int R\t|\pa_j\pa_i\rho|^2dxd\tau
=-\int\div(\r u)\pa_j u_i\pa_j\pa_i\r\,dx\\
\di\quad+\int\pa_j(\r\pa_j u_i)(u\cdot\na\pa_i\rho+\r\pa_i\div u+\pa_i u\cdot\na\r+\pa_i\r\div u)dx\\
\di\quad -\int\Big[u\cdot\na\pa_j u_i+R\pa_j\pa_i\t+\pa_j u\cdot\na u_i+\pa_j\Big(R\f{\t}{\r}\Big)\pa_i\r\Big]\r\pa_j\pa_i\r\,dx\\
\di\quad-\int \f{\pa_j\r}{\r}\pa_j\pa_i\r\big(\mu\Delta u_i+(\mu+\lambda)\pa_i\div u\big) \,dx
+(2\mu+\lambda)\int \pa_j\pa_i\r\pa_j\pa_i\div u\,dx,
\end{array}
\end{equation}
where we have used the following two facts:
\be
\begin{array}{ll}
\di\int \r \pa_j u_{it}\pa_j\pa_i\r \,dx=\f{d}{dt}\int \r\pa_j u_i \pa_j\pa_i\r \,dx
-\int\r_t\pa_j u_i\pa_j\pa_i\r\,dx-\int\r\pa_j u_i\pa_j\pa_i\r_t\,dx\\
\di=\f{d}{dt}\int \r\pa_j u_i \pa_j\pa_i\r \,dx+\int\div(\r u)\pa_j u_i \pa_j\pa_i\r\,dx
+\int \pa_j(\r\pa_j u_i)\pa_i\r_t \,dx\\
\di=\f{d}{dt}\int \r\pa_j u_i \pa_j\pa_i\r\, dx+\int\div(\r u)\pa_j u_i \pa_j\pa_i\r\,dx\\
\di\quad-\int\pa_j(\r\pa_j u_i)(u\cdot\na\pa_i\r+\r\pa_i\div u+\pa_i u\cdot\na\r+\pa_i\r\div u)dx
\end{array}
\ee
and
\be
\di \int\f{1}{\r}\big(\mu\Delta\pa_j u_i+(\mu+\lambda)\pa_j\pa_i\div u\big)\r\pa_j\pa_i\r\,dx
=(2\mu+\lambda)\int\pa_j\pa_i\r\pa_j\pa_i\div u \, dx.
\ee

Next, multiplying the equation \eqref{rho-x1}$_1$ by $\f{2\mu+\lambda}{\r}\pa_j\pa_i\rho$ and integrating with respect to $x$ yield
\be\label{rho-x3}
\ba
 \f{d}{dt}\int\f{2\mu+\lambda}{2\r}|\pa_j\pa_i\r|^2dx
&=-(2\mu+\lambda)\int\pa_j\pa_i\r\pa_j\pa_i\div u \,dx+(2\mu+\lambda)\int\f{\div u}{\r}|\pa_j\pa_i\r|^2dx\\
& -\int \f{2\mu+\lambda}{\r}\pa_j\pa_i\r\big(\pa_j u\cdot\na\pa_i\r+\pa_i u\cdot\na\pa_j\r+\div u\pa_j\pa_i\r\big)dx\\
& -\int \f{2\mu+\lambda}{\r}\pa_j\pa_i\r\big(\pa_j\r\pa_i\div u+\pa_j\pa_i u\cdot\na\r+\pa_i\r\pa_j\div u\big)dx.
\ea
\ee

Finally, we add \eqref{rho-x2} and \eqref{rho-x3} together, sum $i,j$ from $1$ to $3$ and integrate the resulted equation over $(0,t)$ to give
\be\label{rho-x4}
\begin{array}{l}
\di\quad \int\Big(\f{2\mu+\lambda}{2\r}|\na^2\r|^2+\r\pa_j u_i\pa_j\pa_i\r\Big)dx\Big|_{\tau=0}^{\tau=t}
+\int_0^t\int R\t|\na^2\r|^2dxd\tau\\
\di=-\int_0^t\int\div(\r u)\pa_j u_i\pa_j\pa_i\r\, dxd\tau
+\int_0^t\int \pa_j(\r\pa_j u_i)(u\cdot\na\pa_i\r+\r\pa_i\div u+\pa_i u\cdot\na\r\\
\di +\pa_i\r\div u)dxd\tau-\int_0^t\int\Big[u\cdot\na\pa_j u_i+R\pa_j\pa_i\t+\pa_j u\cdot\na u_i+\pa_j\Big(R\f{\t}{\r}\Big)\pa_i\r\Big]\r\pa_j\pa_i\r\,dxd\tau\\
\di-\int_0^t\int \f{\pa_j\r}{\r}\pa_j\pa_i\r\big(\mu\Delta u_i+(\mu+\lambda)\pa_i\div u\big) \,dxd\tau
+(2\mu+\lambda)\int_0^t\int\f{\div u}{\r}|\na^2\r|^2dxd\tau\\
\di -\int_0^t\int \f{2\mu+\lambda}{\r}\pa_j\pa_i\r\big(\pa_j u\cdot\na\pa_i\r+\pa_i u\cdot\na\pa_j\r+\div u\pa_j\pa_i\r\big)dxd\tau\\
\di -\int_0^t\int \f{2\mu+\lambda}{\r}\pa_j\pa_i\r\big(\pa_j\r\pa_i\div u+\pa_j\pa_i u\cdot\na\r+\pa_i\r\pa_j\div u\big)dxd\tau:=\sum_{i=9}^{15}I_i.
\end{array}
\ee
By H\"older's inequality, Cauchy's inequality, Sobolev's inequality, Lemma \ref{appu} and assumption \eqref{assumption}, it holds
\be\label{I9}
\ba
|I_9|&\leq C\int_0^t\int|(\na\r,\na u)||\na u||\na^2\r|dxd\tau\leq C\int_0^t\int(|\na\r|^2+|\na u|^2)|\na^2\r|dxd\tau\\
&\leq C\int_0^t(\|\na\r\|^2_{L^4}+\|\na u\|^2_{L^4})\|\na^2\r\|d\tau
\leq C\int_0^t(\|\na\r\|_{H^1}^2+\|\na u\|_{H^1}^2)\|\na^2\r\|d\tau\\
& \leq \leq C(\chi+\v)\int_0^t(\|(\na^2\r,\na\phi)\|^2+\|(\na^2 u,\na\psi)\|^2)d\tau+C\v.
\ea
\ee
Similar to $I_9$, we have
\be\label{I10}
\di |I^1_{10}|\leq \int_0^t\int(|\na\r||\na u|+|\na^2 u|)|\na^2\r|dxd\tau
\leq \big(\s+C(\chi+\v)\big)\int_0^t\|\na^2\r\|^2d\tau+C_{\s}\int_0^t\|\na^2 u\|^2d\tau.
\ee
It follows from Young's inequality, Sobolev's inequality, Lemma \ref{appu} and assumption \eqref{assumption} that,
\be
\ba
|I^2_{10}|&\leq \int_0^t\int(|\na\r||\na u|+|\na^2 u|)|\na^2 u|dxd\tau
\leq C\int_0^t\big(\|\na\r\|_{L^4}\|\na u\|_{L^4}\|\na^2 u\|+\|\na^2 u\|^2\big)d\tau\\
&\leq C\int_0^t\|\na\r\|_{H^1}\|\na u\|_{H^1}\|\na^2 u\|d\tau+C\int_0^t\|\na^2 u\|^2d\tau\\
&\leq C(\chi+\v)\int_0^t\|(\na^2\r,\na\phi)\|^2d\tau+C\int_0^t\|\na^2 u\|^2d\tau+C\v,
\ea
\ee
\be
\ba
|I^3_{10}|+|I^4_{10}|&\leq C\int_0^t\int(|\na\r||\na u|+|\na^2 u|)|\na\r||\na u|dxd\tau\\
&\leq C\int_0^t\|\na\r\|^2_{L^4}\|\na u\|^2_{L^4}d\tau+C\int_0^t\|\na^2 u\|\|\na\r\|_{L^4}\|\na u\|_{L^4}d\tau\\
&\leq C(\chi+\v)\int_0^t(\|(\na^2\r,\na\phi)\|^2+\|(\na^2 u,\na\psi)\|^2)d\tau+C\v,
\ea
\ee
and
\be\label{I11-1}
|I^1_{11}|+|I^2_{11}|\leq \s\int_0^t\|\na^2\r\|^2d\tau+C_{\s}\int_0^t\|\na^2(u,\t)\|^2d\tau.
\ee
Similar to $I_9$, one has
\be
\ba
|I^3_{11}|+|I^4_{11}|&\leq \int_0^t\int(|\na u|^2+|\na(\r,\t)||\na\r|)|\na^2\r|dxd\tau
\leq C\int_0^t\|\na(\r,u,\t)\|^2_{L^4}\|\na^2\r\|d\tau\\
&\leq C(\chi+\v)\int_0^t(\|\na^2\r\|^2+\|\na^2u\|^2+\|\na^2\t\|^2)d\tau.
\ea
\ee
By H\"older's inequality, Sobolev's inequality, Young's inequality, Lemma \ref{appu} and assumption \eqref{assumption}, we have
\be
\ba
|I_{12}|+|I_{15}|&\leq C\int_0^t\|\na\r\|_{L^4}\|\na^2 u\|_{L^4}\|\na^2\r\|d\tau
\leq C\int_0^t\|\na\r\|_{H^1}\|\na^2\r\|\|\na^2 u\|_{H^1}d\tau\\
&\leq C(\chi+\v)\int_0^t(\|\na^2\r\|^2+\|(\na^3 u,\na^2u)\|^2)d\tau.
\ea
\ee
The same as \eqref{hard}, it holds
\be\label{I13}
|I_{13}|+|I_{14}|\leq C(\chi+\v)\int_0^t(\|\na^2\r\|^2+\|\na^3 u\|^2)d\tau.
\ee
Substituting \eqref{I9}-\eqref{I13} into \eqref{rho-x4} leads to
\be
\begin{array}{l}
\di\|\na^2\r(t)\|^2+\int_0^t\|\na^2\r\|^2d\tau\leq C\|\na u_0,\na^2\r_0\|^2+C\|\na u(t)\|^2\\
\di\quad\quad\quad\quad+C\int_0^t\|\na^2(u,\t)\|^2d\tau+C(\chi+\v)\int_0^t\|\na^3u\|^2d\tau,
\end{array}
\ee
which along with Lemma \ref{appu} and \eqref{psi-1x} implies \eqref{phi-xx}, and the proof of Lemma \ref{le-3}
is completed. \hfill $\Box$

\

Thus, Lemma \ref{high-1} and \ref{le-3} imply

\begin{lemma}\label{le-4}
For $T>0$ and $(\phi,\psi,\z)\in X(0,T)$ satisfying a priori assumption \eqref{assumption} with suitably small $\chi+\v$, we have for $t\in[0,T]$,
\be\label{2nd}
\begin{array}{l}
\di\quad \|(\na\psi,\na\z,\na^2\phi)(t)\|^2+\int_0^t\|\na^2(\phi,\psi,\z)\|^2d\tau\\
\di\leq C(\|(\phi_0,\psi_0,\z_0)\|_{H^1}^2+\|\na^2\phi_0\|^2+\v^{\f18})+C(\chi+\v)\int_0^t\|\na^3 u\|^2d\tau.
\end{array}
\ee
\end{lemma}

\

Finally, we want to derive the highest order derivatives of $\psi$ and $\z$. It holds

\begin{lemma}\label{le-5}
For $T>0$ and $(\phi,\psi,\z)\in X(0,T)$ satisfying a priori assumption \eqref{assumption} with suitably small $\chi+\v$, it holds for $t\in[0,T]$,
\be\label{psi-x1}
\begin{array}{l}
\di\quad \|\na^2(\psi,\z)(t)\|^2+\int_0^t\|\na^3(\psi,\z)\|^2d\tau
\leq C\|(\phi_0,\psi_0,\z_0)\|_{H^2}^2+C\v^{\f18}.
\end{array}
\ee
\end{lemma}

\

{\textbf{Proof}:} First, applying $\pa_i~(i=1,2,3)$ to the equation \eqref{ns-1}$_{2}$ gives
\be
\begin{array}{l}
\di \pa_i u_t+u\cdot\na\pa_i u+R\f{\t}{\r}\na\pa_i\r+R\na\pa_i\t+\pa_i u\cdot\na u+\pa_i\Big(R\f{\t}{\r}\Big)\na\r\\
\di=\f{1}{\r}\big(\mu\Delta\pa_i u+(\mu+\lambda)\na\pa_i\div u\big)-\f{\pa_i\r}{\r^2}\big(\mu\Delta u+(\mu+\lambda)\na\div u\big).
\end{array}
\ee
Multiplying the above equation by $(-\Delta\pa_i u)$, similar to \eqref{psi-2}, we have
\be\label{psi-x2}
\begin{array}{l}
\di\quad \int\f{|\na^2 u|^2}{2}dx\Big|_{\tau=0}^{\tau=t}+\int_0^t\int\Big(\f{\mu}{\r}|\na^3 u|^2+\f{\mu+\lambda}{\r}|\na^2\div u|^2\Big)dxd\tau\\
\di=\int_0^t\int \Big[u\cdot\na\pa_i u+R\f{\t}{\r}\na\pa_i\r+R\na\pa_i\t+\pa_i u\cdot\na u+\pa_i\Big(R\f{\t}{\r}\Big)\na\r\Big]\cdot\Delta\pa_i u \,dxd\tau\\
\di +\int_0^t\int\f{1}{\r^2}\big(\mu\na\r\cdot\na\pa_j\pa_i u\cdot\pa_j\pa_i u-\mu\pa_j\r\pa_j\pa_i u\cdot\Delta\pa_i u
+(\mu+\lambda)\na\r\cdot\na\pa_i\div u\pa_i\div u\\
\di -(\mu+\lambda)\Delta\pa_i u\cdot\na\r \pa_i\div u\big)dxd\tau
+\int_0^t\int\f{\pa_i\r}{\r^2}\big(\mu\Delta u+(\mu+\lambda)\na\div u\big)\cdot\Delta\pa_i u\,dxd\tau:=\sum_{i=16}^{18}I_i.
\end{array}
\ee
It follows from Young's inequality, Sobolev's inequality, Lemma \ref{appu} and assumption \eqref{assumption} that
\be\label{I16}
\di |I_{16}|\leq \big(\s+(\chi+\v)\big)\int_0^t\|\na^3 u\|^2d\tau+C_{\s}\int_0^t\|\na^2(\r,u,\t)\|^2d\tau
\ee
and
\be\label{I17}
\ba
|I_{17}|+|I_{18}|&\leq C\int_0^t\|\na\r\|_{L^4}\|\na^2 u\|_{L^4}\|\na^3u\|d\tau
\leq C\int_0^t\|\na\r\|_{H^1}\|\na^2u\|_{H^1}\|\na^3u\|d\tau\\
&\leq  C(\chi+\v)\int_0^t(\|\na^2 u\|^2+\|\na^3 u\|^2)d\tau.
\ea
\ee
Substituting \eqref{I16}-\eqref{I17} into \eqref{psi-x2} yields
\be\label{psi-x3}
\di \|\na^2 u(t)\|^2+\int_0^t\|\na^3 u\|^2d\tau\leq C\|\na^2 u_0\|^2
+C\int_0^t\|\na^2(\r,u,\t)\|^2d\tau.
\ee

\

Next, applying $\pa_i~(i=1,2,3)$ to the equation \eqref{ns-1}$_3$ gives
\be\label{zeta-x1}
\begin{array}{l}
\di \f{R}{\gamma-1}(\pa_i\t_t+u\cdot\na\pa_i\t)+R\t\pa_i\div u+\f{R}{\gamma-1}\pa_iu\cdot\na\t+R\pa_i\t\div u
=\f{\k}{\r}\Delta\pa_i\t\\
\di+\f{1}{\r}\big(\mu(\na u + (\na u)^t)\cdot\pa_i(\na u + (\na u)^t)+2\lambda\div u\pa_i\div u\big)
-\f{\pa_i\r}{\r^2}\big(\k\Delta\t+\frac{\mu}{2}|\na u + (\na u)^t|^2 + \lambda(\div u)^2\big).
\end{array}
\ee
Multiplying the above equation by $(-\Delta\pa_i\t)$, similar to \eqref{psi-x2}, we have
\be\label{zeta-x2}
\begin{array}{l}
\di \f{R}{\gamma-1}\int\f{|\na^2\t|^2}{2}dx\Big|_{\tau=0}^{\tau=t}+\int_0^t\int\f{\k}{\r}|\na^3\t|^2dxd\tau
=\int_0^t\int\Big[\f{R}{\gamma-1}u\cdot\na\pa_i\t+R\t\pa_i\div u\\
\di\quad+\f{R}{\gamma-1}\pa_iu\cdot\na\t+R\pa_i\t\div u
\Big]\Delta\pa_i\t\,dxd\tau+\int_0^t\int\f{\k}{\r^2}(\na\r\cdot\na\pa_j\pa_i\t\pa_j\pa_i\t-\pa_j\r\pa_j\pa_i\t\Delta\pa_i\t\\
\di\quad+\pa_i\r\Delta\t\Delta\pa_i\t)dxd\tau - \int_0^t\int\f{1}{\r}\big(\mu(\na u + (\na u)^t)\cdot\pa_i(\na u + (\na u)^t)+2\lambda\div u\pa_i\div u\big)\Delta\pa_i\t\,dxd\tau\\
\di\quad +\int_0^t\int\f{\pa_i\r}{\r^2}\big(\frac{\mu}{2}|\na u + (\na u)^t|^2 + \lambda(\div u)^2\big)\Delta\pa_i\t\,dxd\tau:=\sum_{i=19}^{22}I_i.
\end{array}
\ee
By Young's inequality, Sobolev's inequality, Lemma \ref{appu} and assumption \eqref{assumption}, one has
\be\label{I19}
|I_{19}|\leq \big(\s+(\chi+\v)\big)\int_0^t\|\na^3\t\|^2d\tau+C_{\s}\int_0^t\|\na^2(u,\t)\|^2d\tau.
\ee
Similar to \eqref{I17}, it holds,
\be
|I_{20}|\leq C\int_0^t\|\na\r\|_{L^4}\|\na^2 \t\|_{L^4}\|\na^3\t\|d\tau
\leq C(\chi+\v)\int_0^t(\|\na^3\t\|^2+\|\na^2\t\|^2)d\tau.
\ee
It follows from H\"older's inequality, Sobolev's inequality, Young's inequality, Lemma \ref{appu} and assumption \eqref{assumption} that
\be
\ba
|I_{21}|&\leq C\int_0^t\|\na u\|_{L^4}\|\na^2 u\|_{L^4}\|\na^3 \t\|d\tau
\leq C\int_0^t\|\na u\|_{H^1}\|\na^2 u\|_{H^1}\|\na^3 \t\|d\tau\\
&\leq C(\chi+\v)\int_0^t(\|\na^3(u,\t)\|^2+\|\na^2 u\|^2)d\tau,
\ea
\ee
and
\be\label{I22}
\ba
|I_{22}|&\leq C\int_0^t\|\na \r\|_{L^4}\|\na u\|^2_{L^8}\|\na^3 \t\|d\tau
\leq C\int_0^t\|\na \r\|_{H^1}\|\na u\|_{H^2}\|\na^3 \t\|d\tau\\
&\leq C(\chi+\v)\int_0^t(\|\na^3(u,\t)\|^2+\|(\na^2 u,\na\psi)\|^2)d\tau+C\v.
\ea
\ee
Substituting \eqref{I19}-\eqref{I22} into \eqref{zeta-x2} gives
\be\label{zeta-x3}
\di \|\na^2\t(t)\|^2+\int_0^t\|\na^3\t\|^2d\tau\leq C\|\na^2\t_0\|^2+C\int_0^t\|\na^2(u,\t)\|^2d\tau
+C(\chi+\v)\int_0^t\|\na^3 u\|^2d\tau.
\ee

Combining \eqref{psi-x3} and \eqref{zeta-x3}, we derive
\be
\di \|\na^2(u,\t)(t)\|^2+\int_0^t\|\na^3(u,\t)\|^2d\tau\leq C\|\na^2(u_0,\t_0)\|^2
+C\int_0^t\|\na^2(\r,u,\t)\|^2d\tau,
\ee
which along with Lemma \ref{appu} and \eqref{2nd} leads to \eqref{psi-x1}, and the proof of Lemma \ref{le-5} is completed.
\hfill $\Box$

\

{\textbf{Proof of Proposition \ref{priori}}:}  Combining \eqref{basic}, \eqref{phi-x}, \eqref{2nd} and \eqref{psi-x1} together,
we can obtain \eqref{full-es}, the proof of Proposition \ref{priori} is completed.
\hfill $\Box$

\

%
%
%

\section{Proof of Theorem \ref{thm}}\label{main}
\setcounter{equation}{0}

{\textbf{Proof of Theorem \ref{thm} }:}  We now finish the proof of the main result in Theorem \ref{thm}. The global
existence result follows immediately from Proposition \ref{priori} (A priori estimates) and local existence which can be obtained similarly as in \cite{Nash} and \cite{Solo}.
To complete the proof of Theorem \ref{thm}, we only need to justify the time-asymptotic behavior \eqref{large-time}.
In fact, from the estimates \eqref{full-es}, it holds that
\be\label{4.1}
\int_0^{\infty}\Big(\|\na(\phi,\psi,\z)\|^2+\Big|\f{d}{dt}\|\na(\phi,\psi,\z)\|^2\Big|\Big)d\tau<\infty,
\ee
which implies
\be\label{4.2}
\di \lim_{t\rightarrow\infty}\|\na(\phi,\psi,\z)(t)\|^2=0.
\ee
By three-dimensional Sobolev's inequality, one has
\be
\begin{array}{l}
\di \|(\phi,\psi,\z)(t)\|^2_{L^{\infty}}\leq C\|(\phi,\psi,\z)(t)\|\|\na(\phi,\psi,\z)(t)\|+ C\|\na(\phi,\psi,\z)(t)\|\|\na^2(\phi,\psi,\z)(t)\|,
\end{array}
\ee
which together with \eqref{full-es} and \eqref{4.2} yields
$$
\lim_{t\rightarrow\infty}\|(\phi,\psi,\z)(t)\|_{L^{\infty}}=0.
$$
Hence we obtain \eqref{large-time} and finish the proof of Theorem \ref{thm}.

\hfill $\Box$

\noindent\textbf{Acknowledgment.} The research of T. Wang is partially supported by NNSFC grant No. 11601031. The research of Y. Wang is partially supported by the NNSFC grants No. 11671385 and 11688101 and CAS Interdisciplinary Innovation Team.


\begin{thebibliography}{00}

\bibitem{BCK}
J. Brezina, E. Chiodaroli and O. Kreml,
\emph{On contact discontinuities in multi-dimensional isentropic Euler
equations}, preprint, https://arxiv.org/abs/1707.00473.

\bibitem{CC}
G. Q. Chen and J. Chen,
\emph{Stability of rarefaction waves and vacuum states for the multidimensional Euler equations}, J. Hyperbolic Differ. Equ. \textbf{4} (1) (2007), 105-122.

\bibitem{CDK-1}
E. Chiodaroli, C. De Lellis and O. Kreml,
\emph{Global ill-posedness of the isentropic system of gas dynamics},
Comm. Pure Appl. Math.  \textbf{68} (7) (2015), 1157-1190.

\bibitem{CK}
E. Chiodaroli and O. Kreml,
\emph{Non-uniqueness of admissible weak solutions to the Riemann problem for the isentropic Euler equations}, preprint, https://arxiv.org/abs/1704.01747.

\bibitem{DS}
C. De Lellis and L. Sz$\acute{{\rm e}}$kelyhidi Jr.,
\emph{The Euler equations as a differential inclusion},
Ann. of Math.  \textbf{170} (3) (2009), 1417-1436.

\bibitem{FK}
E. Feireisl and O. Kreml,
\emph{Uniqueness of rarefaction waves in multidimensional compressible Euler system},
J. Hyperbolic Differ. Equ. \textbf{12} (3) (2015), 489-499.

\bibitem{FKV}
E. Feireisl, O. Kreml, and  A. Vasseur,
\emph{Stability of the isentropic Riemann solutions of the full multi-dimensional Euler system},
SIAM J. Math. Anal., \textbf{47} (3) (2015), 2416-2425.

\bibitem{goodman}
J. Goodman,
\emph{Nonlinear asymptotic stability of viscous shock profiles for conservation laws},
Arch. Rational. Mech. Anal., \textbf{95} (1986),  325-344.

\bibitem{HM}
H. Hokari and  A. Matsumura,
\emph{Asymptotics toward one-dimensional rarefaction wave for the solution of two-dimensional compressible Euler equation with an artificial viscosity},
Asymptot. Anal. \textbf{15} (1997), 283-298.

\bibitem{Huang-Li-matsumura}
F. M. Huang, J. Li and  A. Matsumura
\emph{ Asymptotic stability of combination of viscous contact wave with rarefaction waves for one-dimensional compressible Navier-Stokes system},
 Arch. Rational Mech. Anal., \textbf{197} (2010), 89-116.

\bibitem{Huang-matsumura}
 F. M. Huang and A. Matsumura,
\emph{ Stability of a composite wave of two viscous shock waves for the full compressible Navier-Stokes equation},
 Comm. Math. Phys., \textbf{289} (2009), 841-861.

\bibitem{Huang-Matsumura-Xin}
 F. M. Huang, A. Matsumura and Z. P. Xin,
 \emph{Stability of contact discontinuities for the 1-D compressible Navier-Stokes equations},
 Arch. Rational Mech. Anal., \textbf{179} (2006), 55-77.

\bibitem{huang-wang}
F. M. Huang and T. Wang,
\emph{Stability of superposition of viscous contact wave and rarefaction waves for compressible Navier-Stokes system},
 Indiana U. Math. J., \textbf{65} (2016), 1833-1875.

%
%
%



\bibitem{Huang-Xin-Yang}
 F. M. Huang, Z. P. Xin and T. Yang,
 \emph{Contact discontinuities with general perturbation for gas motion},
 Adv. Math., \textbf{219} (2008), 1246-1297.


\bibitem{Ito}
K. Ito,
\emph{Asymptotic decay toward the planar rarefaction waves of solutions for viscous conservation laws in several space dimensions},
Math. Models and Methods in Appl. Scis. \textbf{6} (3) (1996), 315-338.


\bibitem{JWX}
Q. S. Jiu, Y. Wang and Z.P. Xin,
\emph{Vacuum behaviors around rarefaction waves to 1D compressible Navier-Stokes equations with density-dependent viscosity},
SIAM J. Math. Anal., \textbf{45} (2013),  3194-3228.

\bibitem{KM}
C. Klingenberg and S. Markfelder,
\emph{The Riemann problem for the multi-dimensional isentropic system of gas dynamics is ill-posed if it contains a shock}, to appear in Arch Rational Mech Anal (2017). https://doi.org/10.1007/s00205-017-1179-z.




\bibitem{Lax}
P. Lax,
\emph{Hyperbolic systems of conservation laws, II},
Comm. Pure Appl. Math., {\bf 10} (1957), 537-566.

\bibitem{LW}
L. A. Li, Y. Wang, Stability of the planar rarefaction wave to the two-dimensional compressible Navier-Stokes equations, Preprint, https://arxiv.org/abs/1710.06063.

\bibitem{LWW}
M. J. Li, T. Wang, Y. Wang, The limit to rarefaction wave with vacuum for 1D compressible fluids with temperature-dependent transport coefficients. Anal. Appl. (Singap.), {\bf 13 } (2015),  555-589.

\bibitem{Liu-1985}
T. P. Liu,
\emph{Nonlinear stability of shock waves for viscous conservation laws},
Mem. Amer. Math. Soc., {\bf 56} (1985), 1-108.

\bibitem{liu-xin-1}
T.P. Liu and Z. P. Xin,
\emph{Nonlinear stability of rarefaction waves for compressible Navier-Stokes equations},
Comm. Math. Phys. \textbf{118} (1988), 451-465.


\bibitem{Liu-Xin}
 T. P. Liu and Z. P.  Xin,
 \emph{Pointwise decay to contact discontinuities for systems of viscous conservation laws},
  Asian J. Math., \textbf{1} (1997), 34-84.

\bibitem{liu-yu-rare}
T. P. Liu and  S. H. Yu,
 \emph{Viscous rarefaction waves},
 Bull. Inst. Math. Acad. Sin. (N.S.) \textbf{5} 2 (2010), 123-179.

\bibitem{Liu-Zeng}
T. P. Liu and Y. N. Zeng,
\emph{Shock waves in conservation laws with physical viscosity},
Mem. Amer. Math. Soc., \textbf{234} (2015), No. 1105.

\bibitem{MN-85}
 A. Matsumura and K. Nishihara,
 \emph{On the stability of traveling wave solutions of a one-dimensional model system for compressible viscous gas},
 Japan J. Appl. Math., \textbf{2}  (1985), 17-25.

\bibitem{MN-86}
 A. Matsumura and K. Nishihara,
 \emph{Asymptotics toward the rarefaction wave of the solutions of a one-dimensional model system for compressible viscous gas},
 Japan J. Appl. Math., \textbf{3} (1986), 1-13.

\bibitem{MN-92}
A. Matsumura and K. Nishihara,
\emph{Global stability of the rarefaction wave of a one-dimensional model system for compressible viscous gas},
Commun. Math. Phys. \textbf{144} (2) (1992), 325-335.

\bibitem{Nash} J. Nash, Le probl\`{e}me de Cauchy pour les \'{e}quations diff\'{e}rentielles d'un fluide g\'{e}n\'{e}ral, \textit{Bull. Soc. Math. France}, \textbf{90}, 487-497 (1962).


\bibitem{NYZ}
K. Nishihara, T. Yang and H. J. Zhao,
\emph{Nonlinear stability of strong rarefaction wave for compressible Navier-Stokes equations},
SIAM J. Math. Anal., \textbf{35} (2004), 1561-1597.

\bibitem{NN}
M. Nishikawa and K. Nishihara,
\emph{Asymptotics toward the planar rarefaction wave for viscous conservation
law in two space dimensions},
Trans. Amer. Math. Soc. \textbf{352} (3) (2000),  1203-1215.


\bibitem{Smoller} 
 J. Smoller,
 \emph{``Shock Waves and Reaction-Diffusion Equations"},
 New York: Springer, 1994.

\bibitem{Szepessy-Xin}
A. Szepessy and Z. P. Xin,
\emph{Nonlinear stability of viscous shock waves},
Arch. Rational Mech. Anal., {\bf 122} (1993), 53-103.

\bibitem{Solo}
V. A. Solonnikov,
\emph{On solvability of an initial-boundary value problem for the equations of motion of a viscous compressible fluid}. in: Studies on Linear Operators and Function Theory.  \textbf{6} [in Russain],
Nauka, Leningrad, (1976), 128-142.

\bibitem{Xin}
Z. P. Xin,
\emph{Asymptotic stability of planar rarefaction waves for viscous conservation laws in several dimensions}, Trans. Amer. Math. Soc. \textbf{319} (1990), 805-820.

\bibitem{xin-2}
Z.P. Xin,
\emph{On nonlinear stability of contact discontinuities. In: Hyperbolic problems:
theory, numerics, applications}  (Stony Brook, NY, 1994), 249-257. World Sci. Publishing,
River Edge, NJ, 1996


\bibitem{zum-de}
K. Zumbrun and D. Serre,
\emph{Viscous and Inviscid Stability of Multidimensional Planar Shock Fronts},
Indiana U. Math. J., \textbf{48} (3) (1999), 937-992.

\end{thebibliography}
\end{document}